\documentclass[10pt]{amsart}

\usepackage[foot]{amsaddr}
\usepackage[utf8]{inputenc} 
\usepackage[T1]{fontenc}
\usepackage{url}
\usepackage[margin=1.3in]{geometry}
\usepackage{parskip}

\usepackage{amssymb} 
\usepackage{amsthm} \usepackage{amsmath}
\usepackage{mathabx} 
\usepackage{mathtools}
\usepackage{mathdots}

\usepackage{relsize} 
\usepackage{exscale} 
\usepackage{tikz}
\usepackage{tikz-cd}
\usetikzlibrary{backgrounds} 
\usetikzlibrary{external}
\usepackage{multirow}

\usepackage{cite}
\makeatletter
\newcommand{\citecomment}[2][]{\citen{#2}#1\citevar}
\newcommand{\citeone}[1]{\citecomment{#1}}
\newcommand{\citetwo}[2][]{\citecomment[,~#1]{#2}}
\newcommand{\citevar}{\@ifnextchar\bgroup{;~\citeone}{\@ifnextchar[{;~\citetwo}{]}}}
\newcommand{\citefirst}{\@ifnextchar\bgroup{\citeone}{\@ifnextchar[{\citetwo}{]}}}
\newcommand{\cites}{[\citefirst}
\makeatother

\usepackage{upref} 

\usepackage[section]{algorithm}
\usepackage{algpseudocode}

\newcommand{\oldalglinenumber}{}
\NewDocumentCommand{\NoLNFor}{ m }{%
  \RenewCommandCopy{\oldalglinenumber}{\alglinenumber}
  \RenewDocumentCommand{\alglinenumber}{ m }{}
  \For{#1}
  \addtocounter{ALG@line}{-1}
  \RenewCommandCopy{\alglinenumber}{\oldalglinenumber}
}
\NewDocumentCommand{\NoLNEndFor}{}{%
  \RenewCommandCopy{\oldalglinenumber}{\alglinenumber}
  \RenewDocumentCommand{\alglinenumber}{ m }{}
  \EndFor
  \addtocounter{ALG@line}{-1}
  \RenewCommandCopy{\alglinenumber}{\oldalglinenumber}
}
\NewDocumentCommand{\NoLNIf}{ m }{%
  \RenewCommandCopy{\oldalglinenumber}{\alglinenumber}
  \RenewDocumentCommand{\alglinenumber}{ m }{}
  \If{#1}
  \addtocounter{ALG@line}{-1}
  \RenewCommandCopy{\alglinenumber}{\oldalglinenumber}
}
\NewDocumentCommand{\NoLNEndIf}{}{%
  \RenewCommandCopy{\oldalglinenumber}{\alglinenumber}
  \RenewDocumentCommand{\alglinenumber}{ m }{}
  \EndIf
  \addtocounter{ALG@line}{-1}
  \RenewCommandCopy{\alglinenumber}{\oldalglinenumber}
}

\usepackage{float}

\allowdisplaybreaks 

\usepackage{enumerate}
\usepackage[shortlabels]{enumitem}

\newcommand{\R}{\mathbb{R}}

\newcommand{\Z}{\mathbb{Z}}

\newcommand{\F}{\mathbb{F}}

\newcommand{\inv}{^{-1}}

\newcommand{\spl}{Spli\textbf{F}f}

\newcommand{\FU}{\mathbb{F}[U, U^{-1}]}
\newcommand{\tatb}{\mathbb{T}_\alpha \cap \mathbb{T}_\beta}
\newcommand{\HFK}{\widehat{HFK}}

\usepackage{graphicx} 
\usepackage{svg} 
\usepackage{tikz}
\usepackage[center]{caption} 

\usepackage[hidelinks,bookmarksopen, bookmarksopenlevel=0]{hyperref}
\usepackage{bookmark}

\newtheorem{theor}{Theorem} 
\newtheorem{thm}{Theorem}[section]
\newtheorem{prop}[thm]{Proposition}

\newtheorem{lemma}[thm]{Lemma}
\newtheorem{cor}[thm]{Corollary}

\newtheorem{conj}{Conjecture}[section]
\newtheorem{claim}{Claim}
\newtheorem*{thm*}{Theorem}
\newtheorem*{prop*}{Proposition}
\newtheorem*{cor*}{Corollary}
\newtheorem*{conj*}{Conjecture}

\newtheorem*{theorintrocfk}{Theorem ~\ref{theor:introcfk}}
\newtheorem*{theorcompute}{Theorem ~\ref{theor:compute}}
\newtheorem*{propth2hfkspliff}{Proposition ~\ref{prop:th2hfkspliff}}

\theoremstyle{definition}
\newtheorem{defn}[thm]{Definition}

\theoremstyle{remark}

\title[Computation of the knot Floer complex of knots of thickness one]{Computation of the knot Floer complex \\of knots of thickness one}
\author{Patricia Sorya}
\address{Département de mathématiques, Université du Québec à Montréal}
\email{sorya.patricia@courrier.uqam.ca}

\thanks{This work was supported by the FRQNT under doctoral research grant 305903}

\date{\today}

\begin{document}

\begin{abstract}
We develop and implement an algorithm that computes the full knot Floer complex of knots of thickness one. As an application, by extending this algorithm to certain knots of thickness two, we show that all but finitely many non-integral Dehn surgery slopes are characterizing for most knots with up to 17 crossings.
\end{abstract}
\maketitle

\newcommand{\perc}{95.83\%}
\newcommand{\percunsolved}{4.17\%}
\newcommand{\thonecomputed}{2 376 438}
\newcommand{\thonespliff}{2 199 957}
\newcommand{\thoneperc}{88.03\%}
\newcommand{\thtwocomputed}{1634}
\newcommand{\thtwospliff}{1489}
\newcommand{\thtworemaining}{49 675}


\section{Introduction}

Knot Floer homology, introduced by Rasmussen \cite{rass} and independently by Ozsváth and Szabó \cite{osHF}, is a knot invariant that has proven to be effective for studying various topological properties of knots in $S^3$, such as fibredness, genus and concordance. It can be obtained from a richer algebraic structure, the \textit{knot Floer complex}. This complex retains more data about the knot, providing further invariants, some of which are particularly useful for the study of Dehn surgeries.

While there are available algorithms for computing knot Floer homology, there is currently no implemented algorithm that effectively outputs the knot Floer complex of an arbitrary knot in $S^3$. The grid diagram algorithm of Manolescu, Ozsváth and Sarkar \cite{mos} has led to a program that calculates knot Floer homology \cite{baldwingillam}, but the high number of generators it considers makes it impractical for the computation of the full knot Floer complexes. Another knot Floer homology calculator, developed by Ozsváth and Szabó \cite{os_uv}, uses bordered algebras to provide more information about the knot Floer complex, but it only yields a quotiented version rather than the full complex.

In this paper, we present and implement an algorithm that recovers the full knot Floer complex of any knot of thickness at most one in $S^3$, from the quotiented complex of Ozsváth and Szabó. 

\begin{theor}\label{theor:introcfk}
The full knot Floer complex of a knot of thickness at most one is determined by the data of its horizontal and vertical arrows.
\end{theor}

The algorithm is grounded in the work of Popovi\'{c} \cite{popovic} who classified the direct sum components of knot Floer complexes of knots of thickness one. The proof of this classification has Theorem \ref{theor:introcfk} as a consequence.

We apply our algorithm to the study of characterizing Dehn surgeries. We show that for the vast majority of knots with up to 17 crossings, all but finitely many non-integral Dehn surgeries are characterizing. This supports McCoy's conjecture asserting the same statement for all knots \cite[Conjecture 1.1]{mcc_hfk}.

\begin{theor}\label{theor:compute}
    Out of the 9 755 329 prime knots with at most 17 crossings, at least \perc{} admit only finitely many non-integral non-characterizing Dehn surgeries.
\end{theor}

This result is achieved by computationally verifying an algebraic condition formulated by \mbox{McCoy}, \textit{property \spl{}}, concerning the homology modules $A^+_k$ of the knot Floer complex. We first identify knots whose knot Floer homology is simple enough to guarantee this condition, by using McCoy's previous work for knots of thickness at most one \cite[Corollary 1.4, Proposition 1.6]{mcc_hfk} and the following proposition for thickness-two knots.

\begin{prop}\label{prop:th2hfkspliff}
Let $K$ be a knot of thickness two. Let $\rho$ be an integer such that for all $s$, the knot Floer homology group $\HFK_d(K,s)$ is non-zero only for gradings $d \in \{s+\rho, s+\rho-1, s+\rho-2\}$.

    Suppose $\rho \in \{0,1,2\}$. If for each $k \geq 0$, at least one of the groups $\HFK_{k+\rho}(K,k)$ or $\HFK_{k+\rho-2}(K,k)$ is trivial, then $K$ and its mirror both satisfy property \spl{}. Therefore, $K$ admits only finitely many non-integral non-characterizing Dehn surgeries.
\end{prop}

We then compute the structure of the modules $A^+_k$ for most of the knots that do not verify \cite[Proposition 1.6]{mcc_hfk} or Proposition \ref{prop:th2hfkspliff}. For thickness-one knots, this is done by using our algorithm to compute the full knot Floer complex, from which we extract the modules $A^+_k$. For thickness-two knots, we adapt the algorithm to recover sufficient information about the modules $A^+_k$ and apply it to cases within our computational capabilities. In particular, for all knots with up to 16 crossings, our strategy yields the full knot Floer complex due to the work of Hanselman who computationally verified, using immersed curves, that the statement of Theorem \ref{theor:introcfk} holds for these knots \cite[Corollary 12.6]{hanselman}. We note that Hanselman's computation also provides a description of their knot Floer complex, as immersed curves turn out to capture the necessary structure for these knots.

Furthermore, our computation showcases the limitations of McCoy's algebraic condition in addressing \cite[Conjecture 1.1]{mcc_hfk}, with the remaining \percunsolved{} of unresolved cases providing examples of knots that do not satisfy property \spl{}. Notably, this includes knots of thickness one, whereas previously identified examples had thickness at least two \cite[Proposition 3.3(ii), Example 3.4]{mcc_hfk}. 

\subsection{Structure of paper} 
The paper is organized as follows. In Section \ref{section:algebraic}, we introduce the algebraic settings in which knot Floer complexes will be studied. Section \ref{section:uniquelift} contains the proof of Theorem \ref{theor:introcfk}. In Section \ref{section:overview}, we present an overview of the algorithm for computing the knot Floer complex of knots of thickness at most one. Section \ref{section:matricial} translates the problem into a computational framework where the differential map is encoded as a matrix. In Section \ref{section:d2=0}, we show that certain degree constraints reduce the problem to a system of linear equations. Section \ref{section:implement} describes the SageMath implementation of the algorithm. In Sections \ref{section:application}, \ref{section:finite1} and \ref{section:finite2} we extend and apply our algorithm to study characterizing Dehn surgeries.

\subsection{Acknowledgements}
I would like to thank David Popovi\'{c}, Jennifer Hom, Ina Petkova and Jonathan Hanselman for interesting discussions, as well as Duncan McCoy and Steven Boyer for their guidance throughout this work.

I extend my gratitude to Franco Saliola for sponsoring my access to the computing platform of Calcul Québec, and to their support staff for excellent assistance. I also thank Cédric Beaulac for tips on structuring the presentation of an algorithm. Lastly, I am deeply grateful to Dan Radulescu for invaluable advice on coding and algorithmic design.

\section{Algebraic setting}\label{section:algebraic}

Knot Floer complexes come in a variety of algebraic flavours. We are interested in the \textit{full knot Floer complex}, from which all other variants can be derived. This full complex can itself be described in different algebraic settings. We present two such settings and we show that the data they encode is equivalent.

\subsection{Basic construction}
We first recall the basics of the construction of a knot Floer complex. From a knot $K$ in $S^3$, we obtain a doubly pointed Heegaard diagram $\mathcal{H} = (\Sigma, \alpha, \beta, w, z)$, where $\Sigma$ is a genus-$g$ surface, $\alpha$ and $\beta$ are sets of $g$ curves on $\Sigma$ and $w, z$ are the two basepoints. A knot Floer complex for $K$ associated to $\mathcal{H}$ is generated by $\tatb = (\alpha_1 \times \ldots \times \alpha_g) \cap (\beta_1 \times \ldots \times \beta_g)$ in the $g$-fold symmetric product $\operatorname{Sym}^g(\Sigma)$. The differential of a knot Floer complex counts certain representatives of Whitney discs $\phi \in \pi_2(x,y)$ between two generators $x,y \in \tatb$, and their intersections with certain auxiliairy submanifolds associated to the basepoints.

In this section, we will assume that any knot Floer complex mentioned refers to a fixed knot $K$ and is obtained from a fixed Heegaard diagram $\mathcal{H}$ for $K$. The knot Floer complex is an invariant of $K$ up to filtered chain homotopy equivalence and does not depend on the choice of $\mathcal{H}$. Therefore, instead of writing $CFK^\infty(\mathcal{H})$ for instance, we may simply write $CFK^\infty(K)$.

\subsection{Knot Floer complex as an \texorpdfstring{$\F[U, U\inv]$}{F[U,U\1]}-module}
We now recall the classical presentation of the knot Floer complex $CFK^\infty(K)$ as an $\F[U, U\inv]$-module from \cite{osHF}, a knot invariant up to filtered homotopy equivalence. $\F$ denotes the field with two elements and $U$ is a formal variable. Let $CFK^-(K)$ be the chain complex generated by $\tatb$ over the ring $\F[U]$ with differential given by
\[d x = \sum_{y \in \tatb}\sum_{\substack{\phi \in \pi_2(x,y) \\ \mu(\phi)=1}} \#(\mathcal{M}(\phi)/\R) \cdot U^{n_w(\phi)} y,\]
where $\mathcal{M}(\phi)$ is the moduli space of holomorphic representatives of the Whitney disc $\phi$, $\mu(\phi)$ is the expected dimension of $\mathcal{M}(\phi)$, and $n_w(\phi)$ is the algebraic intersection number of $\phi$ with $\{w\}\times \operatorname{Sym}^{g-1}(\Sigma)$. The chain complex $CFK^\infty(K)$ is defined as $CFK^-(K) \otimes_{\F[U]}\FU$.

We may visually depict a representative of $CFK^\infty(K)$ in a $\Z \oplus \Z$ lattice as follows. An element $U^i x, x \in \tatb$, has position $(-i, A(U^i x))$, where $A(U^i x)$ is the Alexander grading of $U^i x$. We have in fact $A(U^i x) = A(x)-i$, so all elements $U^i x, i \in \Z$ are represented on a diagonal line of slope 1 intersecting the vertical axis at $A(x)$. If there is a Whitney disc $\phi \in \pi_2(x,y), \mu(\phi)=1$, then $A(U^i x) - A(U^i y) = n_z(\phi) - n_w(\phi)$.

Homogeneous elements of $CFK^\infty(K)$ are endowed with an additional grading called the Maslov grading. The action of multiplication by $U$, modifies this grading by $-2$, i.e. $M(U^i x) = M(x) - 2i$. If there is a Whitney disc $\phi \in \pi_2(x,y), \mu(\phi)=1$, then $M(U^i x) - M(U^i y) = 1 - 2n_w(\phi)$. Thus, the differential lowers the Maslov grading by 1, making it the homological degree on $CFK^\infty(K)$. Therefore, we may interchangeably use \textit{grading} and \textit{degree} to refer to the Maslov grading.

Arrows are drawn between generators to indicate the differential. Arrows are said to be \textit{horizontal}, \textit{vertical} or \textit{diagonal} with respect to this visual representation. The position of an element in the $\Z \oplus \Z$ lattice indicates its filtration level, with respect to the partial order on $\Z \oplus \Z$ given by
\[(i,j) \leq (i',j') \Longleftrightarrow i\leq i' \text{ and } j \leq j',\]
with a strict inequality if $i < i'$ or $j < j'$.

The filtered chain homotopy class of $CFK^\infty(K)$ can be represented by a reduced chain complex (see for instance \cite[Section 2.1]{heddenwatson}). Let $x$ and $y$ be generators of a reduced representative $(C,d)$ of $CFK^\infty(K)$ such that $U^k y$ has non-zero coefficient in $d(U^i x)$. Since the differential $d$ strictly lowers the filtration, we have $U^k y < U^i x$. Therefore, $-k \leq -i$ and $A(y)-k \leq A(x)-i$, where $-k < -i$ or $A(y)-k < A(x)-i$.

\pgfdeclarelayer{edgelayer}
\pgfdeclarelayer{nodelayer}
\pgfsetlayers{background,edgelayer,nodelayer,main}

\tikzstyle{none}=[inner sep=0mm]
\tikzstyle{generator}=[fill=black, draw=none, shape=circle, inner sep=0pt, minimum size=4pt]

\tikzstyle{differential}=[stealth-, line width=0.75pt]
\tikzstyle{axis}=[<-, line width=0.25pt, draw={rgb,255: red,191; green,191; blue,191}]
\tikzstyle{grid}=[-, dotted, line width=0.25pt, draw={rgb,255: red,191; green,191; blue,191}]
\begin{figure}[ht!]
	\centering
	\scalebox{0.75}{%
\begin{tikzpicture}%
	\begin{pgfonlayer}{nodelayer}
		\node [style=generator, label={above:$Ux_0$}] (14) at (-2, 2) {};
		\node [style=generator, label={above left:$Ux_1$}] (15) at (-2, 0) {};
		\node [style=generator, label={below right:$U^{-1}x_5$}] (16) at (2, 0) {};
		\node [style=generator, label={right:$U^{-1}x_3$}] (17) at (2, 2) {};
		\node [style=generator, label={left:$x_2$}] (18) at (0, 2) {};
		\node [style=generator, label={below right:$x_4$}] (19) at (0, -2) {};
		\node [style=generator, label={right:$U^{-1}x_6$}] (20) at (2, -2) {};
		\node [style=none] (21) at (0, 6) {};
		\node [style=none] (22) at (0, -6) {};
		\node [style=none] (23) at (-6, 0) {};
		\node [style=none] (24) at (6, 0) {};
		\node [style=none] (25) at (-1, 6) {};
		\node [style=none] (26) at (-3, 6) {};
		\node [style=none] (27) at (-5, 6) {};
		\node [style=none] (28) at (1, 6) {};
		\node [style=none] (29) at (3, 6) {};
		\node [style=none] (30) at (5, 6) {};
		\node [style=none] (31) at (5, -6) {};
		\node [style=none] (32) at (3, -6) {};
		\node [style=none] (33) at (1, -6) {};
		\node [style=none] (34) at (-1, -6) {};
		\node [style=none] (35) at (-3, -6) {};
		\node [style=none] (36) at (-5, -6) {};
		\node [style=none] (37) at (-6, 5) {};
		\node [style=none] (38) at (-6, 3) {};
		\node [style=none] (39) at (-6, 1) {};
		\node [style=none] (40) at (-6, -1) {};
		\node [style=none] (41) at (-6, -3) {};
		\node [style=none] (42) at (-6, -5) {};
		\node [style=none] (43) at (6, 5) {};
		\node [style=none] (44) at (6, 3) {};
		\node [style=none] (45) at (6, 1) {};
		\node [style=none] (46) at (6, -1) {};
		\node [style=none] (47) at (6, -3) {};
		\node [style=none] (48) at (6, -5) {};
		\node [style=generator, label={above:$x_0$}] (49) at (0.5, 4.5) {};
		\node [style=generator, label={above left:$x_1$}] (50) at (0.5, 2.5) {};
		\node [style=generator, label={right:$U^{-2}x_5$}] (51) at (4.5, 2.5) {};
		\node [style=generator, label={right:$U^{-2}x_3$}] (52) at (4.5, 4.5) {};
		\node [style=generator, label={above:$U^{-1}x_2$}] (53) at (2.5, 4.5) {};
		\node [style=generator, label={below right:$U^{-1}x_4$}] (54) at (2.5, 0.5) {};
		\node [style=generator, label={right:$U^{-2}x_6$}] (55) at (4.5, 0.5) {};
		\node [style=generator, label={left:$U^2x_0$}] (56) at (-4.5, -0.5) {};
		\node [style=generator, label={left:$U^2x_1$}] (57) at (-4.5, -2.5) {};
		\node [style=generator, label={below right:$x_5$}] (58) at (-0.5, -2.5) {};
		\node [style=generator, label={below left:$x_3$}] (59) at (-0.5, -0.5) {};
		\node [style=generator, label={left:$Ux_2$}] (60) at (-2.5, -0.5) {};
		\node [style=generator, label={left:$Ux_4$}] (61) at (-2.5, -4.5) {};
		\node [style=generator, label={right:$x_6$}] (62) at (-0.5, -4.5) {};
		\node [style=none] (63) at (-4, -4) {$\iddots$};
		\node [style=none] (64) at (5.5, 5.5) {$\iddots$};
	\end{pgfonlayer}
	\begin{pgfonlayer}{edgelayer}
		\draw [style=grid] (37.center) to (43.center);
		\draw [style=grid] (44.center) to (38.center);
		\draw [style=grid] (39.center) to (45.center);
		\draw [style=grid] (46.center) to (40.center);
		\draw [style=grid] (41.center) to (47.center);
		\draw [style=grid] (48.center) to (42.center);
		\draw [style=grid] (26.center) to (35.center);
		\draw [style=grid] (27.center) to (36.center);
		\draw [style=grid] (25.center) to (34.center);
		\draw [style=grid] (28.center) to (33.center);
		\draw [style=grid] (29.center) to (32.center);
		\draw [style=grid] (30.center) to (31.center);
		\draw [style=axis] (24.center) to (23.center);
		\draw [style=axis] (21.center) to (22.center);
		\draw [style=differential] (15) to (14);
		\draw [style=differential] (15) to (16);
		\draw [style=differential] (19) to (18);
		\draw [style=differential] (19) to (20);
		\draw [style=differential] (16) to (17);
		\draw [style=differential] (18) to (17);
		\draw [style=differential] (15) to (18);
		\draw [style=differential] (19) to (16);
		\draw [style=differential] (50) to (49);
		\draw [style=differential] (50) to (51);
		\draw [style=differential] (54) to (53);
		\draw [style=differential] (54) to (55);
		\draw [style=differential] (51) to (52);
		\draw [style=differential] (53) to (52);
		\draw [style=differential] (50) to (53);
		\draw [style=differential] (54) to (51);
		\draw [style=differential] (57) to (56);
		\draw [style=differential] (57) to (58);
		\draw [style=differential] (61) to (60);
		\draw [style=differential] (61) to (62);
		\draw [style=differential] (58) to (59);
		\draw [style=differential] (60) to (59);
		\draw [style=differential] (57) to (60);
		\draw [style=differential] (61) to (58);
	\end{pgfonlayer}
\end{tikzpicture}%
}
\caption{The complex $CFK^\infty(K)$ for the $(2,-1)$-cable of the left-handed trefoil}\label{fig:cfk_FU}
\end{figure}

\subsection{Knot Floer complex as an \texorpdfstring{$\F[u,v]$}{F[u,v]}-module}
We also recall the presentation of the knot Floer complex $CFK_{\F[u,v]}(K)$ as an $\F[u,v]$-module, also a knot invariant up to homotopy equivalence, as introduced in \cite{zemke_uv} and summarized in \cite{hom_notes}. As before, $\F$ is the field with two elements and $u,v$ are formal variables. The ring $\F[u,v]$ is bigraded by a $u$-grading $gr_u$ and a $v$-grading $gr_v$ such that $(gr_u(u), gr_v(u)) = (-2, 0)$ and $(gr_u(v), gr_v(v)) = (0, -2)$. The complex $CFK_{\F[u,v]}(K)$ is generated by $\tatb$ over $\F[u,v]$ and the differential is given by
\[d_{\F[u,v]} x = \sum_{y \in \tatb}\sum_{\substack{\phi \in \pi_2(x,y) \\ \mu(\phi)=1}} \#(\mathcal{M}(\phi)/\R) \cdot u^{n_w(\phi)} v^{n_z(\phi)} y.\]
A representative of $CFK_{\F[u,v]}(K)$ admits a decomposition into direct summands $\mathcal{A}_s(K), s \in \Z$, consisting of all $\F$-linear combinations of elements $u^i v^j x, x \in \tatb, (i,j) \in \Z \oplus \Z$, that have $\mathcal{A}$-grading $s \in \Z$, where $\mathcal{A}(u^i v^j x)=(gr_u(u^i v^j x)-gr_v(u^i v^j x))/2$. The $\mathcal{A}$-grading $\mathcal{A}(x)$ of a generator $x \in \tatb$ agrees with its Alexander grading $A(x)$. Since the action of multiplication by $u$ modifies the $u$-grading by $-2$ and the multiplication by $v$ leaves it untouched, and vice versa for the $v$-grading, we have $gr_u(u^i v^j x) = gr_u(x) - 2i$ and $gr_v(u^i v^j x) = gr_v(x) - 2j$. The $u$-grading of $u^i v^j x$ agrees with the Maslov grading of $U^i x \in CFK^\infty(K)$ described above.

In a visual representation of $\mathcal{A}_s(K)$ for some $s \in \Z$, an element $u^i v^j x, x \in \tatb$ has relative position $(-i,-j)$ in the $\Z \oplus \Z$ lattice. Arrows are drawn between generators to indicate the differential. This complex has an implicit filtration given by the powers of $u$ and $v$, since by definition, the differential always increases these powers. This agrees with the partial order on $\Z \oplus \Z$ mentioned above.
We may extend this visual representation to the tensor product $\mathcal{A}_s(K) \otimes_{\F[uv]} \F[uv, (uv)\inv]$, which we denote by $CFK^\infty_{\F[u,v], s}(K)$.

\begin{figure}[ht!]
	\centering
	\scalebox{0.75}{%
\begin{tikzpicture}%
	\begin{pgfonlayer}{nodelayer}
		\node [style=generator, label={left:$uv^{-1}x_0$}] (14) at (-2, 2) {};
		\node [style=generator, label={above left:$ux_1$}] (15) at (-2, 0) {};
		\node [style=generator, label={below right:$u^{-1}x_5$}] (16) at (2, 0) {};
		\node [style=generator, label={right:$u^{-1}x_3$}] (17) at (2, 2) {};
		\node [style=generator, label={below right:$v^{-1}x_2$}] (18) at (0, 2) {};
		\node [style=generator, label={below right:$vx_4$}] (19) at (0, -2) {};
		\node [style=generator, label={right:$u^{-1}vx_6$}] (20) at (2, -2) {};
		\node [style=none] (21) at (0, 6) {};
		\node [style=none] (22) at (0, -6) {};
		\node [style=none] (23) at (-6, 0) {};
		\node [style=none] (24) at (6, 0) {};
		\node [style=none] (25) at (-1, 6) {};
		\node [style=none] (26) at (-3, 6) {};
		\node [style=none] (27) at (-5, 6) {};
		\node [style=none] (28) at (1, 6) {};
		\node [style=none] (29) at (3, 6) {};
		\node [style=none] (30) at (5, 6) {};
		\node [style=none] (31) at (5, -6) {};
		\node [style=none] (32) at (3, -6) {};
		\node [style=none] (33) at (1, -6) {};
		\node [style=none] (34) at (-1, -6) {};
		\node [style=none] (35) at (-3, -6) {};
		\node [style=none] (36) at (-5, -6) {};
		\node [style=none] (37) at (-6, 5) {};
		\node [style=none] (38) at (-6, 3) {};
		\node [style=none] (39) at (-6, 1) {};
		\node [style=none] (40) at (-6, -1) {};
		\node [style=none] (41) at (-6, -3) {};
		\node [style=none] (42) at (-6, -5) {};
		\node [style=none] (43) at (6, 5) {};
		\node [style=none] (44) at (6, 3) {};
		\node [style=none] (45) at (6, 1) {};
		\node [style=none] (46) at (6, -1) {};
		\node [style=none] (47) at (6, -3) {};
		\node [style=none] (48) at (6, -5) {};
		\node [style=generator, label={above left:$v^{-2}x_0$}] (49) at (0.5, 4.5) {};
		\node [style=generator, label={above left:$v^{-1}x_1$}] (50) at (0.5, 2.5) {};
		\node [style=generator, label={right:$u^{-2}v^{-1}x_5$}] (51) at (4.5, 2.5) {};
		\node [style=generator, label={right:$u^{-2}v^{-2}x_3$}] (52) at (4.5, 4.5) {};
		\node [style=generator, label={above:$u^{-1}v^{-2}x_2$}] (53) at (2.5, 4.5) {};
		\node [style=generator, label={below right:$u^{-1}x_4$}] (54) at (2.5, 0.5) {};
		\node [style=generator, label={right:$u^{-2}x_6$}] (55) at (4.5, 0.5) {};
		\node [style=generator, label={left:$u^2x_0$}] (56) at (-4.5, -0.5) {};
		\node [style=generator, label={left:$u^2vx_1$}] (57) at (-4.5, -2.5) {};
		\node [style=generator, label={below right:$vx_5$}] (58) at (-0.5, -2.5) {};
		\node [style=generator, label={below left:$x_3$}] (59) at (-0.5, -0.5) {};
		\node [style=generator, label={left:$ux_2$}] (60) at (-2.5, -0.5) {};
		\node [style=generator, label={left:$uv^2x_4$}] (61) at (-2.5, -4.5) {};
		\node [style=generator, label={right:$v^2x_6$}] (62) at (-0.5, -4.5) {};
		\node [style=none] (63) at (-4, -4) {$\iddots$};
		\node [style=none] (64) at (5.5, 5.5) {$\iddots$};
	\end{pgfonlayer}
	\begin{pgfonlayer}{edgelayer}
		\draw [style=grid] (37.center) to (43.center);
		\draw [style=grid] (44.center) to (38.center);
		\draw [style=grid] (39.center) to (45.center);
		\draw [style=grid] (46.center) to (40.center);
		\draw [style=grid] (41.center) to (47.center);
		\draw [style=grid] (48.center) to (42.center);
		\draw [style=grid] (26.center) to (35.center);
		\draw [style=grid] (27.center) to (36.center);
		\draw [style=grid] (25.center) to (34.center);
		\draw [style=grid] (28.center) to (33.center);
		\draw [style=grid] (29.center) to (32.center);
		\draw [style=grid] (30.center) to (31.center);
		\draw [style=axis] (24.center) to (23.center);
		\draw [style=axis] (21.center) to (22.center);
		\draw [style=differential] (15) to (14);
		\draw [style=differential] (15) to (16);
		\draw [style=differential] (19) to (18);
		\draw [style=differential] (19) to (20);
		\draw [style=differential] (16) to (17);
		\draw [style=differential] (18) to (17);
		\draw [style=differential] (15) to (18);
		\draw [style=differential] (19) to (16);
		\draw [style=differential] (50) to (49);
		\draw [style=differential] (50) to (51);
		\draw [style=differential] (54) to (53);
		\draw [style=differential] (54) to (55);
		\draw [style=differential] (51) to (52);
		\draw [style=differential] (53) to (52);
		\draw [style=differential] (50) to (53);
		\draw [style=differential] (54) to (51);
		\draw [style=differential] (57) to (56);
		\draw [style=differential] (57) to (58);
		\draw [style=differential] (61) to (60);
		\draw [style=differential] (61) to (62);
		\draw [style=differential] (58) to (59);
		\draw [style=differential] (60) to (59);
		\draw [style=differential] (57) to (60);
		\draw [style=differential] (61) to (58);
	\end{pgfonlayer}
\end{tikzpicture}%
}
\caption{The complex $CFK^\infty_{\F[u,v], 0}(K)$ for the $(2,-1)$-cable of the left-handed trefoil}\label{fig:cfk_Fuv}
\end{figure}

\subsection{Equivalence between algebraic settings}
The two algebraic settings contain the same information for a given knot, as given by the next proposition.

\begin{prop}\textup{\cite[Section 1.5]{zemke_conventions}}\label{prop:isomCFKFuv-CFKinfty}
    Let $CFK^\infty_{\F[u,v], s}(\mathcal{H}),\ s \in \Z$ and $CFK^\infty(\mathcal{H})$ be representatives of $CFK^\infty_{\F[u,v], s}(K)$ and $CFK^\infty(K)$ respectively, obtained from the same Heegaard diagram $\mathcal{H}$.
    Then each complex $CFK^\infty_{\F[u,v], s}(\mathcal{H}),\ s \in \Z$, is isomorphic to $CFK^\infty(\mathcal{H})$ by an isomorphism that respects both the filtration up to translation and the $\FU$-module structure, by setting $U=uv$.
\end{prop}
\begin{proof}
    We define a map $\varphi_s: CFK^\infty_{\F[u,v], s}(\mathcal{H}) \rightarrow CFK^\infty(\mathcal{H})$ in the following way.
    
    Let $x \in \tatb$, so that $x$ is a generator of both $CFK_{\F[u,v]}(\mathcal{H})$ and $CFK^\infty(\mathcal{H})$. Set $\varphi_s(u^i v^j x)=U^i x$. We show that $\varphi_s$ is a $U^{\pm 1}$-equivariant filtered chain isomorphism realizing the proposition.
    
\begin{enumerate}[1)]
    \item \textit{$\varphi_s$ is injective:}
    
For a fixed $i$, there is only one possible power $j$ of $v$ such that $u^i v^j x \in CFK^\infty_{\F[u,v], s}(\mathcal{H})$, i.e. $\mathcal{A}(u^i v^j x) = s$. Indeed,
\begin{align*}
    s &= \mathcal{A}(u^i v^j x)\\
    &=(gr_u(u^i v^j x)-gr_v(u^i v^j x))/2 \\
    &= (gr_u(x)-2i-gr_v(x)+2j)/2 \\
    &= (j-i)+A(x)
\end{align*}
implies that $-j = A(x)-i-s$.

    \item \textit{$\varphi_s$ is surjective:}

An element $U^i x \in CFK^\infty(\mathcal{H})$ has antecedent $u^i v^{i-A(x)} x \in CFK^\infty_{\F[u,v], s}(\mathcal{H})$.

    \item \textit{$\varphi_s$ preserves the filtration up to translation:}
    
The element $u^i v^{s+i-A(x)}x \in CFK^\infty_{\F[u,v], s}(\mathcal{H})$ and its image $U^i x \in CFK^\infty(\mathcal{H})$ have respective filtration levels $(-i, A(x)-i-s)$ and $(-i, A(x)-i)$. Therefore, $\varphi_s$ translates filtration levels by $(0,s)$.

    \item \textit{$\varphi_s$ is $U^{\pm 1}$-equivariant:}
    
We have $\varphi_s((uv)^{\pm 1} \cdot u^{i} v^{j} x) = \varphi_s(u^{i\pm1} v^{j\pm1} x) = U^{i\pm1}x = U^{\pm 1} \cdot U^i x =  U^{\pm 1} \varphi_s(u^{i} v^{j} x)$.

    \item \textit{$\varphi_s$ is a chain map:}
    
By definition of the differentials, we have
\begin{align*}
    \varphi_s(d_{\F[u,v]} u^{i} v^{j} x) &= \varphi_s \bigg(\sum_{y \in \tatb}\sum_{\substack{\phi \in \pi_2(x,y) \\ \mu(\phi)=1}} \#(\mathcal{M}(\phi)/\R) \cdot u^{i+n_w(\phi)} v^{j+n_z(\phi)} y \bigg) \\
    &= \sum_{y \in \tatb}\sum_{\substack{\phi \in \pi_2(x,y) \\ \mu(\phi)=1}} \#(\mathcal{M}(\phi)/\R) \cdot \varphi_s 
 (u^{i+n_w(\phi)} v^{j+n_z(\phi)} y) \\
    &= \sum_{y \in \tatb}\sum_{\substack{\phi \in \pi_2(x,y) \\ \mu(\phi)=1}} \#(\mathcal{M}(\phi)/\R) \cdot U^{i+n_w(\phi)} y\\
    &= d(U^i x) \\
    &= d_{\F[u,v]} \varphi_s(u^{i} v^{j} x) \qedhere
\end{align*}
\end{enumerate}
\end{proof}

    Applying the reduction lemma of \cite[Section 2.1]{heddenwatson} in a mirrored way to $CFK^\infty_{\F[u,v], s}(\mathcal{H})$ and $CFK^\infty(\mathcal{H})$, we obtain from Proposition \ref{prop:isomCFKFuv-CFKinfty} a $U^{\pm 1}$-equivariant filtered chain isomorphism $\varphi_s: \mathcal{C}_s \rightarrow \mathcal{C}$ between reduced representatives of $CFK^\infty_{\F[u,v], s}(K)$ and $CFK^\infty(K)$ respectively. We can recover the (reduced) summands $\mathcal{A}_s(K)$ of $CFK_{\F[u,v]}$ from $CFK^\infty(K)$ by restricting $\varphi_s\inv$ to elements with filtration $(i,j) \leq (s,0)$.

\subsection{Thickness}
Both algebraic settings contain the data of the \textit{knot Floer homology} $\widehat{HFK}(K)$ of the knot: on one hand, $\widehat{HFK}(K) \cong CFK_{\F[u,v]}(K)/(u,v)$ and on the other hand,  $\widehat{HFK}(K)$ is the homology of the associated graded complex of $CFK^-(K)/U$. Reduced representatives of $CFK_{\F[u,v]}$ and $CFK^\infty(K)$ have generating sets that are in bijection with the generating set of $\widehat{HFK}(K)$. We denote by $\widehat{HFK}(K,a)$ the knot Floer homology of $K$ in Alexander grading $a$.

The thickness of a knot $K$ is defined from $\widehat{HFK}(K) = \oplus_{a \in \Z} \widehat{HFK}(K,a)$.
\begin{defn}
	The \emph{thickness} of a knot $K$ is the number
	\[th(K) = \max\{|(M(x)-A(x)) - (M(y)-A(y))\,|\, x,y \text{ generators of } \widehat{HFK}(K)\}\]
\end{defn}

A low thickness imposes constraints on the possible arrows representing the differential map. We will apply these constraints in the next section, where we focus on knots of thickness one.
\section{Chain homotopy equivalence of lifts}\label{section:uniquelift}

\subsection{Horizontal and vertical arrows}
The algorithm of Ozsváth and Szabó mentioned in the introduction and implemented in the program SnapPy \cite{SnapPy} outputs the quotient of a reduced representative of $CFK_{\F[u,v]}(K)$ by $uv$, for any knot $K$ given as input. From now on, we will assume that all chain complexes mentioned are reduced. In this subsection, we recall how the horizontal and vertical arrows of the full complex $CFK_{\F[u,v]}(K)$ are captured by this quotiented complex for any knot $K$.

\begin{prop}\label{prop:horizvert}
Let $(C,d)$ be a reduced representative of $CFK_{\F[u,v]}(K)$. Then $(C,d)/(uv)$ is obtained from the data of the horizontal and vertical arrows of $(C,d)$. Conversely,  the data of the horizontal and vertical arrows of $(C,d)$ is contained in $(C,d)/(uv)$. 
\end{prop}
\begin{proof}
Let $x$ be a generator of $C$. The differential of $[x]$ in $(C,d)/(uv)$ is given by
\[ [dx] = \sum_{\substack{y \text{ generator} \\ \text{of } C}} c_y \cdot [u^{i_y} v^{j_y} y] \]
for some $c_y \in \F$ and $i_y, j_y \geq 0 \in \Z$. If $i_y$ and $j_y$ are both non-zero, then $[u^{i_y} v^{j_y} y] = 0$ and therefore
\[ [dx] = \sum_{\substack{y \text{ generator} \\ \text{of } C}} \bigg(\sum_{j_y = 0} c_y [u^{i_y}  y]  + \sum_{i_y = 0} c_y [v^{j_y}  y]\bigg), \]
which is precisely the data of horizontal and vertical arrows leaving $x$ in $(C,d)$.

Since $C$ is generated over $\F[u,v]$, this also gives the data of horizontal and vertical arrows leaving $u^i v^j x$ for all $i,j \in \Z$.
\end{proof}

Note that, due to the isomorphism from the discussion following Proposition \ref{prop:isomCFKFuv-CFKinfty}, the arrows of the quotient complex $(C,d)/(uv)$ also provide the data of the horizontal and vertical arrows of a reduced representative of $CFK^\infty(K)$.

To recover the full knot Floer complex from $(C,d)/(uv)$, we need to find a lift of $(C,d)/(uv)$ to a chain complex $(C',d')$ over $\F[u,v]$ which is chain homotopy equivalent to $(C,d)$. By Proposition \ref{prop:horizvert}, if this lift $(C',d')$ is reduced, we know that it must contain the same data of horizontal and vertical arrows as $(C,d)$.

\subsection{Chain homotopy equivalence}
The lifts of $(C,d)/(uv)$ to complexes over $\F[u,v]$, for representatives $(C,d)$ of $CFK_{\F[u,v]}(K)$, may belong to distinct chain homotopy classes. However, when the thickness of $K$ is at most one, all such complexes are in fact equivalent. 

\begin{thm}\label{theor:homotopyequivalence}
    Let $K$ be a knot of thickness at most one and let $(C,d)$ be a reduced representative of $CFK_{\F[u,v]}(K)$. Then all lifts of $(C,d)/(uv)$ to a reduced complex over $\F[u,v]$ are isomorphic.
\end{thm}

Theorem \ref{theor:homotopyequivalence} combined with Proposition \ref{prop:horizvert} immediately implies Theorem \ref{theor:introcfk}.

\begin{theorintrocfk}
The full knot Floer complex of a knot of thickness at most one is determined by the data of its horizontal and vertical arrows. \qed
\end{theorintrocfk}

The case of thickness zero in Theorem \ref{theor:homotopyequivalence} is trivial since all representatives of $CFK_{\F[u,v]}(K)$ contain only horizontal and vertical arrows. In particular, Petkova showed that the chain homotopy class $CFK_{\F[u,v]}(K)$ of a knot of thickness zero is determined by the knot's Alexander polynomial and $\tau$ invariant \cite[Theorem 4]{petkova}. For knots of thickness one, this is a consequence of the proof of the following result of Popovi\'{c}.

\begin{thm}\textup{\cite[Theorem 1.1]{popovic}}\label{theor:popovic}
    Let $K$ be a knot of thickness one. Then $CFK_{\F[u,v]}(K)$ splits uniquely as a direct sum of an $\F[u,v]$-standard complex of thickness at most 1 and trivial local systems, each of which belongs to a specific set of systems $\mathcal{L}$.
\end{thm}

The $\F[u,v]$-standard complex of Theorem \ref{theor:popovic} is an $\F[u,v]$-realization of a standard complex as originally defined in \cite[Definition 4.3]{dhst}. The exact description of the local systems in $\mathcal{L}$ can be found in the statement of \cite[Theorem 1.1]{popovic}, but the key property of $\mathcal{L}$ relevant to our purposes is the following.

\begin{prop}\textup{\cite[Proposition 4.11]{popovic}}\label{prop:popoviclift}
    Let $C$ be a chain complex over $\F[u,v]$ of thickness one and let $L \in \mathcal{L}$ be a local system such that $C/(uv) \cong L/(uv) \oplus A/(uv)$ for some $\F[u,v]$-chain complex $A$. Then $C \cong L \oplus A$.
\end{prop}

\begin{proof}[Proof of Theorem \ref{theor:homotopyequivalence}]
Let $K$ be a knot of thickness one. 
Let $(C,d)$ be a representative of the chain homotopy class of $CFK_{\F[u,v]}$ and let $(C',d')$ be a lift of $(C,d)/(uv)$ over $\F[u,v]$ at the level of chain complexes. Note that $C'=C$ as bigraded $\F[u,v]$-modules, so we may write $(C,d') = (C',d')$.

We decompose the differential map $d$ into $d = H+V+D$, where $H, V$ and $D$ are respectively the horizontal, vertical and diagonal arrows of $d$. Let $d_{uv}=H+V$.
Similarly, we write $d'= H'+V'+D'$ and $d'_{uv}=H'+V'$.

Since $(C,d)/(uv) \cong (C,d')/(uv)$, we have $d_{uv}=d'_{uv}$ by Proposition \ref{prop:horizvert}.

The splitting of Theorem \ref{theor:popovic} is realized by a change of basis $P$ such that $(C, PdP^{-1})$ is a direct sum as in the statement of Theorem \ref{theor:popovic}  (see proofs of \cite[Lemmas 4.12, 4.13 and 4.14]{popovic}). Restricting $d$ to $d_{uv}$, we have that $(C, Pd_{uv}P^{-1})$ is a direct sum of a standard complex of thickness at most one and local systems from $\mathcal{L}$ with the diagonal arrows removed.

Thus, by Proposition \ref{prop:horizvert}, both $(C,PdP\inv)/(uv)$ and $(C,Pd'P\inv)/(uv)$ are isomorphic to the same direct sum $L_1/(uv) \oplus \ldots \oplus L_k/(uv) \oplus S/(uv), k \geq 0$, where $L_i \in \mathcal{L}, i=1, \ldots, k$, and $S$ is an $\F[u,v]$-standard complex of thickness at most one.

By \cite[Algorithm 3.12]{popovicalgo}, the quotiented standard complex $S/(uv)$ has a unique lift over $\F[u,v]$. Applying Proposition \ref{prop:popoviclift} inductively on the number $k$ of local systems in the direct sum, we obtain the isomorphism $(C,PdP\inv) \cong (C,Pd'P\inv)$. Performing the change of basis $P\inv$ yields the isomorphism $(C,d) \cong (C,d')$ as desired.
\end{proof}

\section{Finding a lift: an overview}\label{section:overview}

In this section, we give an overview of our method to find a lift of $CFK_{\F[u,v]}(K)/(uv)$ for knots of thickness one, which we will detail in the following two sections. By Theorem \ref{theor:homotopyequivalence}, this leads to an algorithm that determines the full knot Floer complex of knots of thickness one.  For computational reasons, we pass to the setting of $CFK^\infty(K)$ over the ring $\FU$, for which we only need to consider a single formal variable $U$.

\renewcommand{\algorithmicrequire}{\textbf{Input:}}
\renewcommand{\algorithmicensure}{\textbf{Output:}}

\begin{algorithm}
\caption{$CFK^\infty(K)$ for knots of thickness $\leq 1$}
\label{algo:cfk1}
    \begin{algorithmic}
\Require Knot $K$ with $th(K) \leq 1$
\Ensure Filtered homotopy representative of $CFK^\infty(K)$
    \end{algorithmic}
\end{algorithm}

The main goal is to construct a chain complex $\mathcal{C}=(C,d)$ over $\FU$ such that $\varphi_s\inv(\mathcal{C})/(uv) \simeq CFK^\infty_{\F[u,v],s}(K)/(uv)$ for all $s \in \Z$. Here $\varphi_s\inv$ are the isomorphisms from the discussion following Proposition \ref{prop:isomCFKFuv-CFKinfty}. We say that such a complex $\mathcal{C}$ is a \textit{lift} of $CFK_{\F[u,v]}(K)/(uv)$. 

Since the $\FU$-module $C$ and the vertical and horizontal arrows of the differential $d$ are known from Ozsváth and Szabó's algorithm, we only need to find the diagonal arrows of $d$.

The first main step is to encode the differential map as a matrix. We construct a matrix $d_{var}$ that contains the data of the known vertical and horizontal arrows, along with entries consisting of unknown variables for possible diagonal arrows, considering constraints given by the Alexander and Maslov gradings. This step does not depend on the thickness of the knot and is described in Section \ref{section:matricial}.

The second main step is to determine a value in $\FU$ for each unknown variable in the matrix $d_{var}$ such that the condition $d_{var}^2 = 0$ of a chain complex is satisfied. We thus rewrite $d_{var}^2=0$ as a set of equations to be solved. By construction of $d_{var}$, a solution to these equations will yield a chain complex $\mathcal{C}=(C,d)$ that respects the filtration and degree constraints expected for a knot Floer complex. The complex $\mathcal{C}$ also has the same data of horizontal and vertical arrows as a reduced representative of $CFK_{\F[u,v]}(K)/(uv)$, making it a lift of $CFK_{\F[u,v]}(K)/(uv)$. A key point for the computational feasibility of our algorithm is that, for thickness one knots, the equations coming from $d_{var}^2=0$ are always linear. This is demonstrated in Section \ref{section:d2=0}. A solution is then obtained by basic linear algebra, giving the desired lift of $CFK_{\F[u,v]}(K)/(uv)$.

\section{Matricial representation}\label{section:matricial}

Our first main step is to encode the differential map $d$ as a matrix with placeholders for the unknown entries. The $\FU$-module underlying $CFK^\infty(K)$ is generated by the generators $x_0, x_1, \ldots, x_{n-1}$ of $\widehat{HFK}(K)$ over $\FU$. Thus, $d$ can be represented by an $n \times n$ matrix with values in $\FU$: the $(i,j)$ entry of this matrix is the coefficient $a_{i,j}$ in $d(x_j) = \sum_{i = 0}^{n-1} a_{i,j} x_i$. In fact, since $d$ respects the filtration, all entries $a_{i,j}$ take values in $\F[U]$. From now on, we will denote both the differential and its matrix by $d$.

\subsection{Entries for horizontal and vertical arrows}
We decompose $d$ into $d = H+V+D$, where $H, V$ and $D$ are respectively the horizontal, vertical and diagonal arrows of the differential.

We recover the matrix $H+V$ using the output from Ozsváth and Szabó's algorithm for computing $CFK_{\F[u,v]}(K)/(uv)$. It provides us with the generators $x_0, x_1, \ldots, x_{n-1}$ and their Maslov and Alexander gradings, and tells us if $H+V$ has an arrow from $x_j$ to $U^k x_i$ for some power $k \geq 0$. Since the differential lowers the Maslov grading by 1 and multiplication by $U$ lowers the Maslov grading by 2, we have
\[M(U^k x_i) = M(x_i) - 2k = M(x_j) - 1.\]
Therefore, if Ozsváth and Szabó's algorithm indicates that there is an arrow from a $x_j$ to $U^k x_i$ for some power $k \geq 0$, we set the $(i,j)$ entry of the matrix $H+V$ to be
\[a_{i,j} = U^{(M(x_i) - M(x_j) + 1)/2}.\]

\subsection{Entries for possible diagonal arrows}
Next, we find pairs of generators of $CFK^\infty(K)$ that may be connected by a diagonal arrow. We consider how a differential map affects the Maslov and Alexander gradings. 

A diagonal arrow from $x_j$ to $U^k x_i$ for some power $k \geq 1$ must meet the conditions $M(U^k x_i) - M(x_j) = -1$ and $A(x_i) - A(x_j) < k$.
Thus, for every $(i,j)$ such that 
\begin{itemize}
    \item[(D1)] $(M(x_i) - M(x_j) + 1)/2 \geq 1$ and
    \item[(D2)] $(M(x_i) - M(x_j) + 1)/2 > A(x_i) - A(x_j)$,
\end{itemize}
there could be a diagonal arrow from $x_j$ to $U^k a_{i,j}$, where 
\begin{itemize}
    \item[(D3)] $k = (M(x_i) - M(x_j) + 1)/2$.
\end{itemize}

We construct a placeholder matrix $D_{var}$ in the following way. If $(i,j)$ satisfies both (D1) and (D2), then the $(i,j)$ entry of $D_{var}$ is $U^k a_{i,j}$, where $a_{i,j}$ is an unknown variable with values in $\F$ and $k$ is as in (D3). Otherwise, the entry is zero. We then form the matrix $d_{var} = H+V+D_{var}$ with entries in $\F[U][\{a_{i,j} \;|\; (i,j) \text{ verify (D1) and (D2)}\}]$. We now want to find the values of $a_{i,j}$ for which $d_{var}$ is a differential map for the $\FU$-module underlying $CFK^\infty(K)$.

\section{Solving for \texorpdfstring{$d^2=0$}{d2=0}}\label{section:d2=0}

Setting $d_{var}^2 = 0$, we obtain equations $[d_{var}^2]_{k,l} = 0$, for each $(k,l) \in \{0, \ldots, n-1\}^2$, where the variables $a_{i,j}$ are the unknowns. Finding these solutions is in general computationally challenging as the equations may involve degree-two polynomials in the ring $\F[\{a_{i,j} \;|\; (i,j) \text{ verify (D1) and (D2)}\}]$, with a number of variables $a_{i,j}$ that can be quite large. However, it turns out that for knots of thickness one, the system $[d_{var}^2]_{k,l} = 0$ consists only of linear equations, which can be solved easily with basic linear algebra.

\subsection{Consecutive diagonal arrows} While the methods of Section \ref{section:matricial} can be applied to any knot, we now restrict our study to knots with low thickness to obtain further constraints on the possible diagonal arrows. The goal of this subsection is to show that given certain degree conditions on $\widehat{HFK}(K)$, there cannot be consecutive diagonal arrows in a reduced chain complex representing $CFK^\infty(K)$.

\begin{prop}\label{prop:Dvar2=0}
    Suppose $K$ is a knot of thickness at most two such that $\widehat{HFK}(K, a)$ is supported in at most 2 degrees for all $a \in \Z$. Then $d_{var}=H+V+D_{var}$ as constructed above is such that $D_{var}^2 = 0$.
\end{prop}

Note that, by the definition of thickness, knots of thickness at most one verify the condition of Proposition \ref{prop:Dvar2=0}. Although Algorithm \ref{algo:cfk1} focuses on this case only, the more general statement of Proposition \ref{prop:Dvar2=0} will be applied in later sections.

Under the condition that the thickness is at most two, we obtain the next three lemmas concerning the Alexander and Maslov gradings of generators connected by a diagonal arrow. We will then use the condition on the support of $\widehat{HFK}(K)$ to prove Proposition \ref{prop:Dvar2=0}.

\begin{lemma}\label{lemma:A}
    Suppose $K$ is a knot of thickness at most two and let $U^k a_{i,j}$ be a non-zero entry in $D_{var}$. Then $|A(x_i) - A(x_j)| \leq 1$.
\end{lemma}
\begin{proof}
    Suppose $A(x_i) - A(x_j) \geq 2$. Then (D3) and (D2) yield
    \begin{align*}
         M(x_i) - A(x_i)& =     M(x_j) + 2k-1 - A(x_i) \\
                    & \geq  M(x_j) + 2(A(x_i)-A(x_j)+1)-1 - A(x_i) \\
                    & \geq     M(x_j) - A(x_j) + 3
    \end{align*}
    which implies that $K$ has thickness at least three.

    Suppose $A(x_j) - A(x_i) \geq 2$. Similarly to the argument above, we obtain
    \begin{align*}
     M(x_i) - A(x_i)& =     (M(x_j) + 2k-1) - A(x_i) \\
                & \geq  M(x_j) + 1 - A(x_j) + 2 \\
                & \geq     M(x_j) - A(x_j) + 3. \qedhere
    \end{align*} 
\end{proof}

\begin{lemma}\label{lemma:nw}
    Suppose $K$ is a knot of thickness at most two and let $U^k a_{i,j}$ be a non-zero entry in $D_{var}$. Let $\eta = A(x_i)-A(x_j)$. Then $\eta \in \{-1,0,1\}$ by Lemma \ref{lemma:A} and $k = 1$ when $\eta = -1$ or $0$, and $k = 2$ when $\eta = 1$.
\end{lemma}
\begin{proof}
We have  $M(x_i) - A(x_i) = M(x_j) + 2k-1 - A(x_j) - \eta$, which implies that $2k-1- \eta \leq 2$, hence $k \leq (3+\eta)/2$. 
Since $k \geq 1$, replacing the value of $\eta$ with $-1,0$ or $1$ in $k \leq (3+\eta)/2$ gives the result.
\end{proof}

\begin{lemma}\label{lemma:MA}
    Suppose $K$ is a knot of thickness at most two and let $U^k a_{i,j}$ be a non-zero entry in $D_{var}$. Then $M(x_i)-A(x_i) = M(x_j)-A(x_j)+2$ when $\eta = 1$ or $-1$, and $M(x_i)-A(x_i) = M(x_j)-A(x_j)+1$ when $\eta = 0$.
\end{lemma}
\begin{proof}
Replace $k$ and $\eta$ in $M(x_i)-A(x_i) = (M(x_j)+2k+1)-(A(x_j)+\eta)$ by the pairs given by Lemma \ref{lemma:nw}.
\end{proof}

\begin{proof}[Proof of Proposition \ref{prop:Dvar2=0}]
    Suppose that $D_{var}^2 \neq 0$. This means that there are non-zero entries $U^{k_1} a_{j,k}$ and $U^{k_2} a_{i,j}$ in $D_{var}$ that contribute $U^{k_1+k_2} a_{i,j}a_{j,k}$ to a non-zero entry of $D_{var}^2 \neq 0$.

If $A(x_k) \neq A(x_j)$ or $A(x_j) \neq A(x_i)$, then by Lemma \ref{lemma:MA},  $(A(x_k) - M(x_k)) - (A(x_i)-M(x_i)) = (A(x_k) - M(x_k)) - (A(x_j)-M(x_j)) + (A(x_j) - M(x_j)) - (A(x_i)-M(x_i)) \geq 3$, which contradicts the thickness of $K$ being at most 2.

If $A(x_k) = A(x_j) = A(x_i)$, then by (D3) we have $M(x_i) = M(x_j) + 2k_2 -1 = M(x_j) + 1$ and $M(x_j) = M(x_k) + 2k_1 -1 = M(x_k) + 1$. Hence, the knot Floer homology of $K$ in Alexander grading $A(x_k) = A(x_j) = A(x_i)$ is supported in at least 3 distinct degrees, a contradiction.
\end{proof}

\subsection{Linear system of equations}\label{section:linearsystem}
We now return to the setting of $d_{var}$ and translate the problem of finding lifts of $CFK_{\F[u,v]}(K)/(uv)$ into a system of linear equations.

\begin{prop}\label{prop:linear}
    Suppose $K$ is a knot of thickness at most two such that $\widehat{HFK}(K, a)$ is supported in at most 2 degrees for all $a \in \Z$. Then the entries of $d_{var}^2$ are polynomials of degree at most one in the variables $a_{i,j}$ over $\F[U]$.
\end{prop}
\begin{proof}
By Proposition \ref{prop:Dvar2=0},
    \begin{align*}
        d_{var}^2 &= (H+V+D_{var})^2 \\
            &= (H+V)^2 + (H+V)D_{var} + D_{var}(H+V) + D_{var}^2\\
            &= (H+V)^2 + (H+V)D_{var} + D_{var}(H+V).
    \end{align*}
The result follows from the fact that $(H+V)$ has entries in $\F[U]$, for which the variables $a_{i,j}$ have degree zero, and $D_{var}$ has entries of the form $U^k a_{i,j}$, where the variables $a_{i,j}$ have degree one.
\end{proof}

We may view the entries of $d_{var}^2$ as polynomials in $U$ with coefficients in $\F\langle\{a_{i,j}\}\rangle$. By setting $d_{var}^2=0$, we must have that each coefficient $\sum a_{i_k,j_k}$ of a power of $U$ is equal to zero. We thus obtain a linear system of equations $E = \{\sum a_{i_k,j_k}=0\}$ over $\F$ where the variables $a_{i,j}$ are the unknowns. This system can be represented by a matrix equation $Aa=b$ where $a$ is the vector of variables $a_{i,j}$ to solve for. 

Given a solution $a=a_0$, we replace its values into the corresponding entries of $D_{var}$ to obtain a matrix $D_0=D_{var}(a_0)$. We then build the differential complex $\mathcal{C}_0 = (C, d_0 = H+V+D_0)$, where $C \cong \widehat{HFK}(K) \otimes \FU$ is the $\FU$-module underlying $CFK^\infty(K)$. By Theorem \ref{theor:homotopyequivalence} and Proposition \ref{prop:isomCFKFuv-CFKinfty}, the complex $\mathcal{C}_0$ is a representative of $CFK^\infty(K)$ if $K$ has thickness at most one.
\section{Implementation}\label{section:implement}

The previous discussion has been implemented in SageMath, utilizing SnapPy \cite{SnapPy} as an imported package. SnapPy is used to input the data of a knot, via its integrated census or a planar diagram, and for calling upon the method \texttt{knot\_floer\_homology}, an implementation of Ozsváth and Szabó's algorithm,  to obtain the data of $CFK_{\F[u,v]}/(uv)$.

SageMath can generate polynomial rings and handle symbolic computations over them. This allows us to extract the equations to be solved over the ring $\F[U]$, as described in Section \ref{section:linearsystem}, and to translate them into a matrix equation $Aa=b$ over $\F$.

To obtain a solution to the matrix equation $Aa=b$, we use SageMath's matrix equation solver \texttt{solve\_right} which implements Gaussian elimination over $\F$.
 
\begin{minipage}{\linewidth}
\rule[-3.5pt]{\textwidth}{0.7pt}
\textbf{Algorithm 4.1 }$CFK^\infty(K)$ for knots of thickness $\leq 1$\\
\rule[8.5pt]{\textwidth}{0.4pt}
\vspace{-2em}
\begin{algorithmic}[1]
\Require Knot $K$ with $th(K) \leq 1$
\Ensure Filtered homotopy representative of $CFK^\infty(K)$

\State Obtain $CFK_{\F[u,v]}/(uv)$ and $\widehat{HFK}$ via the \texttt{knot\_floer\_homology(complex=True)} method 

\State Let $\{x_0, \ldots, x_{n-1}\}$ be the generators of $\widehat{HFK}$
\State Generate the matrix $H+V \in M_n(\F[U])$ of horizontal and vertical arrows from $CFK_{\F[u,v]}/(uv)$
\State Initiate a zero $n \times n$ matrix $D_{var}$ and populate it
\NoLNFor{$i,j \in \{0, \dots, n-1\}$}
    \NoLNIf{$(i,j)$ satisfies (D1) and (D2)} set $[D_{var}]_{i,j} = U^k a_{i,j}$, where $k$ is as in (D3)
    \NoLNEndIf
\NoLNEndFor
\State Generate the matrix equation $Aa=b$
\begin{itemize}
    \item Obtain a set of expressions $E$ from the $\F\langle\{a_{i,j}\}\rangle$ coefficients of non-zero entries of the matrix $(H+V+D_{var})^2$
    \item Let $A$ be the matrix with each row consisting of the $\F$ coefficients of the $a_{i,j}$ for an entry in $E$
    \item Let $b$ be the vector of constant terms for each element in $E$
    \item Let $a$ be the vector of unknown variables $a_{i,j}$
\end{itemize}
\State Find a solution $a_0$ via \texttt{solve\_right}
\State Get a matrix $D_0=D_{var}(a_0)$
\State Construct a chain complex $\mathcal{C}_0=(C, d_0=H+V+D_0)$, where $C \cong \widehat{HFK} \otimes \FU$.
\State \Return $\mathcal{C}_0$
\end{algorithmic}\nopagebreak
\rule[8.5pt]{\textwidth}{0.4pt}
\end{minipage}
\section{Finiteness of non-integral non-characterizing slopes: an overview}\label{section:application}

As an application of Algorithm \ref{algo:cfk1}, we investigate the set of characterizing slopes for knots in $S^3$. A Dehn surgery slope is said to be characterizing for a knot $K$ if the orientation-preserving homeomorphism type of its $p/q$-Dehn surgery $S^3_K(p/q)$ determines $K$ up to isotopy. That is, if there is some knot $K'$ such that $S^3_{K'}(p/q) \cong S^3_K(p/q)$ via an orientation-preserving homeomorphism, then $K' = K$. Baker and Motegi asked whether a non-integral slope $p/q$ is characterizing for a hyperbolic knot when $|p|+|q|$ is sufficiently large \cite[Question 5.6]{bakermotegi}. This naturally leads to the question of whether the same holds for any knot in $S^3$.

\begin{conj}\textup{\cite[Conjecture 1.1]{mcc_hfk}}\label{conj:mcc}
   Let $K$ be a knot in $S^3$. Then all but finitely many non-integral slopes are characterizing for $K$.
\end{conj}

Conjecture \ref{conj:mcc} has been shown to hold for thickness-zero knots, L-space knots \cite[Corollary 1.4]{mcc_hfk} and composite knots \cite[Theorem 2]{sorya}. In this paper, we restrict our attention to prime knots of thickness one and two, and show the conjecture to be true for the vast majority of prime knots with at most 17 crossings.

\begin{theorcompute}
    Out of the 9 755 329 prime knots with at most 17 crossings, at least \perc{} admit only finitely many non-integral non-characterizing Dehn surgeries.
\end{theorcompute}

\subsection{Property \spl{}}
A key result towards Theorem \ref{theor:compute} is a sufficient condition on the knot Floer complex $CFK^\infty(K)$ formulated by McCoy, which guarantees that the conjecture holds for a given knot $K$. Let $C_{\{i \geq 0 \vee j \geq k\}}$ be the quotient complex of $CFK^\infty(K)$ represented by homogenous elements with $\Z \oplus \Z$ filtration satisfying $i \geq 0$ or $j \geq k$, and denote its homology by $A^+_k$. Let $\F_d$ denote an $\F$ summand supported in grading $d$.

\begin{defn}{\cite[Definition 1.5]{mcc_hfk}} A knot $K$ has \textit{property \spl{}} if for all $k \in \Z$, the graded $\F[U]$-module $A^+_k$ admits a direct sum decomposition of the form
\begin{equation}\label{eq:spliff}
A' \oplus \F^{n_1}_{d_1} \oplus \F^{n_2}_{d_2},   
\end{equation}
where $n_1, n_1 \geq 0$, $d_1$ is odd, $d_2$ is even and the $\F[U]$-module $A'$ does not contain a summand whose elements are all killed by the $U$-action.
\end{defn}

\begin{thm}\textup{\cite[Theorem 1.2, Theorem 1.3]{mcc_hfk}}\label{theor:mccspliff}
   Let $K$ be a knot in $S^3$ such that both $K$ and its mirror have property \spl{}. Then all but finitely many non-integral slopes are characterizing for $K$.
\end{thm}

Recall that $A^+_k$ admits a decomposition $A^+_k \cong \mathcal{T}_{-2V_k} \oplus A^{red}_k$ for some integer $V_k \geq 0$, where 
$\mathcal{T}_d = \FU/U\mathbb{F}[U]$ and $1$ has even grading $d$. Since $\mathcal{T}_{-2V_k}$ contains elements that are not killed by $U$, and since there is an even grading shift $A^{red}_k[-2k] \cong A^{red}_{-k}$, showing that $A^{red}_k$ decomposes as in (\ref{eq:spliff}) for all $k \geq 0$ is equivalent to saying that $K$ has property \spl{}. We therefore say that the modules $A^+_k$ and $A^{red}_k$ have \textit{property \spl{}} if $A^{red}_k$ admits a decomposition as in (1).

Theorem \ref{theor:compute} is thus obtained by computing the complexes $A^{red}_k, k\geq 0,$ for knots and their mirrors, and verifying whether they satisfy property \spl{}.

\subsection{Summary of results}
\subsubsection{Thickness-one knots}

We applied Algorithm \ref{algo:cfk1} to all knots obtained from SnapPy's\linebreak \texttt{NonalternatingKnotExteriors} iterator for prime knots with up to 16 crossings and most knots in Regina's database \cite{burton} of prime knots with 17 crossings.

Combining the output of Algorithm \ref{algo:cfk1} and McCoy's work on the structure of the modules $A^+_k$ of thickness-one knots \cite[Section 3.3]{mcc_hfk}, we determine whether property \spl{} is satisfied for each of the 437 982 prime thickness-one knots with at most 16 crossings and their mirrors, and for \thonecomputed{} of the 2 516 641 prime thickness-one knots with 17 crossings and their mirrors. We found that \thonespliff{} pairs of such knots and their mirrors have property \spl{}, thus verifying the conjecture for \thoneperc{} of prime thickness-one knots with at most 17 crossings. In particular, Conjecture \ref{conj:mcc} is solved for all prime knots up to 11 crossings, and all but 6 prime knots with 12 crossings, listed in Table \ref{table:12cross} along with their $A^{red}_k$ module which fails to have property \spl{}.

\begin{table}[ht]
    \centering
    \begin{tabular}{|lcc|} \hline
        Knot & $k$ & $A^{red}_k$ \\ \hline
    12n67   & 0 & $\F_0 \oplus \F_2^2$ \\
    m12n89  & 0 & $\F_0 \oplus \F_2^2$ \\
    m12n134 & 0 & $\F_0 \oplus \F_2^2$ \\
    m12n229 & 0 & $\F_0 \oplus \F_2^2$ \\
    m12n244 & 1 & $\F_2 \oplus \F_4$ \\
    m12n639 & 0 & $\F_0 \oplus \F_2^2$ \\\hline
    \end{tabular}
        \caption{Knots with 12 crossings for which Conjecture \ref{conj:mcc} remains unresolved}
    \label{table:12cross}
\end{table}

\subsubsection{Thickness-two knots}
We also extended the strategy of Algorithm \ref{algo:cfk1} to thicker knots and check whether property \spl{} is satisfied for certain knots of thickness two. To do this, we first establish thickness-two analogues of McCoy's results on the structure of the modules $A^+_k$. In particular, Proposition \ref{prop:th2hfkspliff} gives a condition on the knot Floer homology $\HFK(K)$ of a thickness-two knot $K$ that guarantees that it has property \spl{}.

We then apply the extended algorithm to all thickness-two knots with up to 16 crossings and certain thickness-two knots with 17 crossings. Table \ref{table:th2spliff} provides a breakdown of the number of thickness-two knots up to 16 crossings according to whether both the knot and its mirror satisfy property \spl{}, or whether at least one of them does not.

\begin{table}[ht]
    \centering
    \begin{tabular}{|c|cc|}
    \hline
    Crossings &  $K$ and $mK$ \spl{} & $K$ or $mK$ non-\spl{} \\ \hline
    13    & 3 & 0 \\
    14   & 32 & 9 \\
    15   & 256 & 193 \\
    16   & 2058 & 2578 \\\hline
    \end{tabular}
        \caption{Thickness-two knots up to 16 crossings and property \spl{}}
    \label{table:th2spliff}
\end{table}

For knots with 17 crossings, \thtwospliff{} of the \thtwocomputed{} thickness-two knots for which we were able to compute the structure of the modules $A^{+}_k$ verified property \spl{}. The large number of complexes to generate prevented us from carrying out the computation for the remaining \thtworemaining{} thickness-two knots with 17 crossings. This computational limitation, along with the empirical observation that the proportion of knots satisfying property \spl{} decreases as the number of crossings increases, suggest that another strategy must be considered to solve Conjecture \ref{conj:mcc} for an arbitrary knot.

Combining all this with the fact that thickness-zero knots always have property \spl{} \cite[Proposition 1.6]{mcc_hfk} and that all but 7 knots with at most 17 crossings have thickness at most two, we obtain the computational result stated as Theorem \ref{theor:compute}.

\subsection{Organization towards Theorem \ref{theor:compute}} Sections \ref{section:finite1} and \ref{section:finite2} detail the theoretical results and computational methods required to establish Theorem \ref{theor:compute}. Their content is organized as follows. We first explain our strategy to compute $A^{red}_k$ for knots of thickness one in Section \ref{section:finite1}. We then develop the case of thickness-two knots in Section \ref{section:finite2}. We analyze the structure of their modules $A^+_k$ in Subsection \ref{subsec:th2Ak}, and in Subsection \ref{subsec:th2spliff}, we prove Proposition \ref{prop:th2hfkspliff}. In Subsection \ref{subsec:th2compute}, we explain how Algorithm \ref{algo:cfk1} was extended to compute the modules $A^{red}_k$ for certain thickness-two knots, and thus obtain the statement of Theorem \ref{theor:compute}.

\section{Finiteness of non-integral non-characterizing slopes: thickness one}\label{section:finite1}

\subsection{Computing \texorpdfstring{$A^+_k$}{A+k}}

Recall that $\mathcal{T}_d = \FU/U\mathbb{F}[U]$ where $1$ has even grading $d$. Let $\mathcal{T}_{d^-}|_{\leq d^+}$ denote the sub-$\F[U]$-module of $\mathcal{T}_{d^-}$ generated by $U^{(d^- -d^+ +\epsilon)/2}$, i.e.  
\[\mathcal{T}_{d^-}|_{\leq d^+}= \F\langle1, U^{-1}, U^{-2}, \ldots, U^{(d^- -d^+ +\epsilon)/2}\rangle,\]
where $\epsilon = 0$ if $d^+$ is even and $1$ if $d^+$ is odd. If $d^+ < d^-$, then $\mathcal{T}_{d^-}|_{\leq d^+}=0$. Otherwise, the element $1$ has degree $d^-$ and $U^{(d^- -d^+ +\epsilon)/2}$ has degree $d^+ - \epsilon$. In other words, it is the truncation of the tower $\mathcal{T} = \FU$ with lowest degree $d^-$ and highest degree $d^+ - \epsilon$.

Our main object of interest, the $\F[U]$-module $A^+_k$, is the homology group of the complex $C_{\{i \geq 0 \vee j \geq k\}}$, represented by homogeneous elements of $CFK^\infty(K)$ whose $\Z \oplus \Z$ filtration level $(i,j)$  satisfies $i \geq  0$ or $j \geq k$. This complex has has infinitely many generators over the finite field $\F$, which makes it unpractical for computational manipulation. To address this, we consider instead the quotient complex $C_{\{i < 0 \wedge j \geq k\}}$, represented by homogenous elements of $CFK^\infty(K)$ whose $\Z \oplus \Z$ filtration level $(i,j)$  satisfies $i < 0$ and $j \geq k$. This has finitely many generators over $\F$, so it is well suited for computational encoding. Its homology is related to $A^+_k$ by an $\F[U]$-module isomorphism (see the proof of \cite[Lemma 29]{gainullin} or \cite[Lemma 3.2(i)]{nizhang})
\begin{equation}\label{eq:gainullin}
    H_*(C_{\{i < 0 \wedge j \geq k\}}) \cong \mathcal{T}_{-2{V_k}}|_{\leq -2} \oplus A^{red}_k.
\end{equation}

We do not know \textit{a priori} which components of the $\F[U]$-module $H_*(C_{\{i < 0 \wedge j \geq k\}})$ are mapped to $A^{red}_k$ under this (non-canonical) isomorphism.

The following structural lemma will allow us to recover enough information about $A^+_k$ from $H_*(C_{\{i < 0 \wedge j \geq k\}})$ to conclude whether $A^+_k$ has property \spl{}.

\begin{lemma}\textup{\cite[Lemma 3.14]{mcc_hfk}}\label{lemma:structth1}
    Let $K$ be a knot of thickness one. Let $\rho$ be an integer such that for all $s$, the group $\HFK_d(K,s)$ is non-zero only for gradings $d \in \{s+\rho, s+\rho-1\}$. Then for all $k \geq 0$, there exist integers $a,b \geq  0$ such that $A^+_k$ takes the following form
 \[A^+_k = \mathcal{T}_{\min(0,k+\rho-1\pm \epsilon)} \oplus \mathcal{T}_{2k}|_{\leq k+\rho-2 \pm \eta} \oplus \F_{k+\rho-1}^a \oplus \F_{k+\rho-2}^b,\]
 where $\epsilon = 0$ if $k+\rho-1$ is even and $1$ otherwise, and $\eta = 0$ if $k+\rho-2$ is even and $1$ otherwise.
\end{lemma}

\begin{cor}\label{cor:HCth1Spliff}Let $K$ be a knot of thickness one and $\rho, \epsilon, \eta$ be as in Lemma \ref{lemma:structth1}. Then
\[H_*(C_{\{i < 0 \wedge j \geq k\}}) \cong \mathcal{T}_{\min(0,k+\rho-1\pm \epsilon)}|_{\leq -2} \oplus \mathcal{T}_{2k}|_{\leq k+\rho-2 \pm \eta} \oplus \F_{k+\rho-1}^a \oplus \F_{k+\rho-2}^b \qedhere\]
for all $k \in \Z$ and $K$ has property \spl{} if and only if $H_*(C_{\{i < 0 \wedge j \geq \rho-3\}})$  has property \spl{}.
\end{cor}
\begin{proof}
The isomorphism is a direct consequence of combining (\ref{eq:gainullin}) and Lemma \ref{lemma:structth1}. 

Next, we observe that if the even number among $k+\rho-1$ and $k+\rho-2$ is greater than zero, then a component of $H_*(C_{\{i < 0 \wedge j \geq k\}})$ is mapped by (\ref{eq:gainullin}) into $A^{red}_k$ unless it is supported in negative even degrees.

By \cite[Lemma 3.15]{mcc_hfk}, $K$ may fail to have property \spl{} only if $\rho \geq 3$ and $A^+_{\rho-3}$ does not have property \spl{}. In this case, $k+\rho-1 = 2\rho-4$ and thus $k+\rho-2 = 2\rho-5$ are always greater than zero, so $\mathcal{T}_{\min(0,k+\rho-1\pm \epsilon)}|_{\leq -2}=\mathcal{T}_0|_{\leq -2}$ is trivial. Therefore, $H_*(C_{\{i < 0 \wedge j \geq \rho-3\}}) \cong A^{red}_{\rho-3}$, and $K$ has property \spl{} if and only if $H_*(C_{\{i < 0 \wedge j \geq \rho-3\}})$ has property \spl{}.
\end{proof}

\subsection{Implementation in SageMath}\label{section:Akimpl}
The complex $C_{\{i < 0 \wedge j \geq k\}}$ is generated in the following way. Recall that Algorithm \ref{algo:cfk1} outputs a matrix for the differential of $CFK^\infty(K)$ in the basis given by the generating set of $\HFK(K)$. The basis for $C_{\{i < 0 \wedge j \geq k\}}$ is given by 
\[B = \{U^{-i}x \,|\, A(x)+i \geq k, i < 0, x \text{ a generator of } \HFK(K) \}.\]
We index the elements of $B$ by $b_0, \ldots, b_{m-1}$. An element $b_l = U^{-i}x$ is implemented as an object with attributes recording the index $l \in \{0, \ldots, m-1\}$, the power $-i$ of $U$ and the generator $x \in \HFK(K)$.

We then construct the matrix $d \in M_m(\F)$ of the differential of $C_{\{i < 0 \wedge j \geq k\}}$ in this basis, according to the output of Algorithm \ref{algo:cfk1}. To obtain the homology group $H_*(C_{\{i < 0 \wedge j \geq k\}})$, we use SageMath's built-in \texttt{kernel} and \texttt{image} methods. Next, we use SageMath's \texttt{basis} and \texttt{lift} methods to obtain representatives of the basis elements of $H_*(C_{\{i < 0 \wedge j \geq k\}})$ in the coordinates $b_0, \ldots, b_{m-1}$. We then extract the Maslov index of the (homogeneous) element $ \sum_{i \in I} U^{k_i} x_i$ corresponding to a representative $\sum_{j\in J} b_j$ via the associated object parameters.

Finally, to check for property \spl{} according to Corollary \ref{cor:HCth1Spliff}, we need to understand the $\F[U]$-module structure of $H_*(C_{\{i < 0 \wedge j \geq \rho-3\}})$. The latter may fail to have property \spl{} only if there are elements in both gradings $2\rho-4$ and $2\rho-6$. In this situation, we consider a subset $B' \subset B$ consisting of a representative for each element in grading $2\rho-4$. We have that $H_*(C_{\{i < 0 \wedge j \geq \rho-3\}})$, and thus $K$, has property \spl{} if and only if $UB'$ is not entirely contained in the image of $d$. This condition is verified by iterating through the elements $b \in B'$, stopping if $Ub$ is not in the image of $d$.
\section{Finiteness of non-integral non-characterizing slopes: thickness two}\label{section:finite2}

\subsection{Structure of \texorpdfstring{$A^+_k$}{Ak} for thickness-two knots}\label{subsec:th2Ak}

The aim of this section is to describe the general algebraic structure of the modules $A^+_k$ for knots of thickness two by establishing the following analogue of Lemma \ref{lemma:structth1}.

\begin{lemma}\label{lemma:structth2}
    Let $K$ be a knot of thickness two. Let $\rho$ be an integer such that for all $s$, the group $\HFK_d(K,s)$ is non-zero only for gradings $d \in \{s+\rho, s+\rho-1, s+\rho-2\}$. Then for all $k \geq 0$, there exist integers $r,a,b,c \geq  0$ such that $A^+_k$ takes the following form as an $\F[U]$-module
 \[A^+_k = \mathcal{T}_{\min(0,k+\rho-\eta\pm 1)} \oplus \mathcal{T}_{2k}|_{\leq k+\rho-2 \pm 1} \oplus \left(\frac{\F_{k+\rho-1}[U]}{U^2}\right)^r \oplus\ \F_{k+\rho-1}^a \oplus \F_{k+\rho-2}^b \oplus \F_{k+\rho-3}^c,\]
where $\eta = 1 + (k+\rho \mod 2)$, and $\F_{d}[U]$ is the ring $\F[U]$ where $1$ has degree $d$.
\end{lemma}

We first review a basic algebraic fact about homomorphisms of graded $\F[U]$-modules.

\begin{lemma}\label{lemma:basicFU}
Let $f$ be a graded $\F[U]$-module homomorphism with domain $\mathcal{T}_d$, where $d \in 2\Z \cup \{-\infty\}$ and $\mathcal{T}_{-\infty} = \mathcal{T}$. If $f$ is not trivial, then the image of $f$ is isomorphic to $\mathcal{T}_{d'}$, where $d' \geq d$.
\end{lemma}
\begin{proof}
A graded sub-$\F[U]$-module of $\mathcal{T}_d$ is either $\mathcal{T}_d$, the trivial module, or the module generated by an element of $\mathcal{T}_d$ of degree at least $d$. Since the kernel of $f$ is a sub-$\F[U]$-module of $\mathcal{T}_d$, the statement follows from the first isomorphism theorem.
\end{proof}

\begin{proof}[Proof of Lemma \ref{lemma:structth2}]
    The proof is modelled on the proofs of \cite[Lemma 3.14]{mcc_hfk} and \cite[Theorem 1.4]{os_alt}.

    Denote by $C_{\{\mathcal{F}\}}$ the quotient of $CFK^\infty(K)$ represented by homogenous elements whose $\Z \oplus \Z$ filtration levels $(i,j)$ satisfy the constraint $\mathcal{F}$.

    For a $\Z$-graded $\F$-module $M = \bigoplus_{s \in \Z}M_s$, let $M|_{\geq k}=\bigoplus_{s \geq k}M_s$, $M|_{\leq k}=\bigoplus_{s \leq k}M_s$ and $M|_{k}=M_k$.


We have a short exact sequence of complexes over the ring $\F[U]$
\[0 \to C_{\{i \leq -1 \wedge j \leq k - 1\}} \to CFK^\infty(K) \to C_{\{i \geq 0 \vee j \geq k\}} \to 0\]
which induces an exact triangle of $\F[U]$-modules in homology
\begin{equation}\label{triangle:injection}
\begin{tikzcd}
	{H_*(CFK^\infty(K)) \cong \mathcal{T}} && {H_*(C_{\{i \geq 0 \vee j \geq k\}})=A^+_k} \\
	& {H_*(C_{\{i \leq -1 \wedge j \leq k - 1\}})}
	\arrow["j_*",from=1-1, to=1-3]
	\arrow[from=1-3, to=2-2]
	\arrow[from=2-2, to=1-1]
\end{tikzcd}
\end{equation}

Similarly, we have a short exact sequence of complexes over the ring $\F[U]$
\[0 \to C_{\{i \geq 0 \vee j \geq k\}} \to C_{\{i \geq 0\}} \oplus C_{\{j \geq k\}} \to C_{\{i \geq 0 \wedge j \geq k\}} \to 0
\]

which induces an exact triangle of $\F[U]$-modules in homology
\begin{equation}\label{triangle:surjection}
\begin{tikzcd}
	{H_*(C_{\{i \geq 0 \vee j \geq k\}})=A^+_k} && {H_*(C_{\{i \geq 0\}}) \oplus H_*(C_{\{j \geq k\}})\cong\mathcal{T}_0 \oplus \mathcal{T}_{2k}} \\
	& {H_*(C_{\{j \geq k\}})}
	\arrow["i_*",from=1-1, to=1-3]
	\arrow[from=1-3, to=2-2]
	\arrow[from=2-2, to=1-1]
\end{tikzcd}
\end{equation}

We now examine the gradings of the underlying $\F$-modules in (\ref{triangle:injection}) and (\ref{triangle:surjection}) to recover information about the structure of $A^{red}_k$.

\begin{claim}\label{claim:Aredmaxsupport}
$A^{red}_k$ is supported in degrees at most $k+\rho-1$.
\end{claim}

\begin{proof}[Proof of Claim \ref{claim:Aredmaxsupport}]

By definition of $\rho$, homogeneous elements of $C_{\{i \leq -1 \wedge j \leq k - 1\}}$ have degree at most $k+\rho-2$. Therefore, $H_s(C_{\{i \leq -1 \wedge j \leq k - 1\}}) = 0$ for all $s \geq k+\rho-1$. Forgetting about the $\F[U]$-module structure, the exact triangle (\ref{triangle:injection}) gives rise to a long exact sequence of $\F$-modules with isomorphisms
\[H_{s+1}(CFK^\infty(K)) \cong H_{s+1}(C_{\{i \geq 0 \vee j \geq k\}})\]
for all $s+1 \geq k+\rho$. Thus, we obtain a commutative diagram of graded $\F$-modules
\begin{equation}\label{diag:injectioniso}
\begin{tikzcd}
H_*(CFK^\infty(K))|_{\geq k+\rho} \arrow[d, "\cong"'] \arrow[r, "\cong"] & H_*(C_{\{i \geq 0 \vee j \geq k\}})|_{\geq k+\rho} \arrow[d, "="] \\
\mathcal{T}|_{\geq k+\rho} \arrow[r, "\cong", "j_*|"']                                                      & A^+_k|_{\geq k+\rho}         
\end{tikzcd}
\end{equation}

Recall that $A^+_k$ decomposes as $\mathcal{T}_{-2V_k} \oplus A^{red}_k$. By Lemma \ref{lemma:basicFU}, $j_*(\mathcal{T})=\mathcal{T}_{-2V_k}$.
The isomorphism $j_*|$ of (\ref{diag:injectioniso}) implies that the elements of $A^+_k$ of degree at least $k+\rho$ are precisely those of the $\mathcal{T}_{-2V_k}$ summand. Therefore, the homogeneous elements of $A^{red}_k$ have degree at most $k+\rho-1$. 
\end{proof}


\begin{claim}\label{claim:AredT2k}
All elements of $A^{red}_k$ of degrees at most $k+\rho-4$ belong to a (possibly trivial) $\F[U]$-module summand of the form $\mathcal{T}_{2k}|_{\leq k+\rho-2 \pm 1}$.
\end{claim}
\begin{proof}[Proof of Claim \ref{claim:AredT2k}]
By definition of $\rho$, homogeneous elements of $C_{\{i \geq 0 \wedge j \geq k\}}$ have degree at least $k+\rho-2$. Therefore, $H_s(C_{\{i \geq 0 \wedge j \geq k\}}) = 0$ for all $s \leq k+\rho-3$. Forgetting about the $\F[U]$-module structure, the exact triangle (\ref{triangle:surjection}) becomes a long exact sequence of $\F$-modules which gives isomorphisms
\[H_{s-1}(C_{\{i \geq 0 \vee j \geq k\}}) \cong H_{s-1}(C_{\{i \geq 0\}}) \oplus H_{s-1}(C_{\{j \geq k\}})\]
for all $s-1 \leq k+\rho-4$. Thus, we obtain a commutative diagram of $\F$-modules
\begin{equation}\label{diag:surjectioniso}
    \begin{tikzcd}
H_*(C_{\{i \geq 0 \vee j \geq k\}})|_{\leq k+\rho-4} \arrow[d, "="'] \arrow[r, "\cong"] & (H_*(C_{\{i \geq 0\}}) \oplus H_*(C_{\{j \geq k\}}))|_{\leq k+\rho-4} \arrow[d, "\cong"] \\
A^+_k|_{\leq k+\rho-4} \arrow[r, "\cong", "i_*|"']                                                      & \mathcal{T}_{0}|_{\leq k+\rho-4} \oplus \mathcal{T}_{2k}|_{\leq k+\rho-4}         
\end{tikzcd}
\end{equation}

In grading $k+\rho-3$, the long exact sequence of $\F$-modules underlying (\ref{triangle:surjection}) gives a surjection
\begin{equation}\label{diag:surjection}
\begin{tikzcd}
A^+_k|_{k+\rho-3} \arrow[r, two heads] & \mathcal{T}_{0}|_{k+\rho-3} \oplus \mathcal{T}_{2k}|_{k+\rho-3}
\end{tikzcd}.
\end{equation}

By Lemma \ref{lemma:basicFU}, $i_*(\mathcal{T}_{-2V_k})$ is isomorphic to a sub-$\F[U]$-module $\mathcal{T}_{d}$ of $\mathcal{T}_{0} \oplus \mathcal{T}_{2k}$. Hence, $i_*(\mathcal{T}_{-2V_k})$ is isomorphic to either $\mathcal{T}_{0}$ or $\mathcal{T}_{2k}$

If $k+\rho-3 \geq 0$ then $\mathcal{T}_{-2V_k}|_{\leq k+\rho-3} \neq 0$. The isomorphism $i_*|$ of (\ref{diag:surjectioniso}) tells us that if $i_*(\mathcal{T}_{-2V_k}) \cong \mathcal{T}_{2k}$, then $2k=-2V_k$. Since $2k \geq 0$ and $-2V_k\leq 0$, we have $2k=0$. Therefore, up to composition with an isomorphism, we may assume that  $i_*(\mathcal{T}_{-2V_k}) = \mathcal{T}_{0}$ and $(i_*)\inv(\mathcal{T}_{2k}) \subset A^{red}_k$. This means that the elements of degrees at most $k+\rho-4$ in $A^{red}_k$ are precisely those of a (possibly trivial) truncated tower $\mathcal{T}_{2k}|_{\leq l}$ for some $l \geq k+\rho-4$. Considering the surjection of (\ref{diag:surjection}), we have in fact $l \geq k+\rho-3$.
Claim \ref{claim:Aredmaxsupport} implies that $l \leq k+\rho-1$, from which we conclude that $l=k+\rho-2 \pm 1$.

If $k+\rho-3 < 0$, then $\mathcal{T}_{0}|_{\leq k+\rho-3} \oplus \mathcal{T}_{2k}|_{\leq k+\rho-3} = 0$. By the isomorphism $i_*|$, there are no elements in gradings less than or equal to $k+\rho-4$ in $A^+_k$. We may therefore take $l = k+\rho-3$ so that $\mathcal{T}_{2k}|_{\leq l} = 0$.
\end{proof}

From Claims \ref{claim:Aredmaxsupport} and \ref{claim:AredT2k}, we deduce that $A^{red}_k/(\mathcal{T}_{2k}|_{\leq l})$ is supported in degrees between $k+\rho-3$ and $k+\rho-1$. Therefore, aside from $\mathcal{T}_{2k}|_{\leq l}$, the remaining $\F[U]$-summands of $A^{red}_k$ may only take the forms
\[\frac{\F_{k+\rho-1}[U]}{U^2} ,\  \F_{k+\rho-1} ,\ \F_{k+\rho-2} ,\ \F_{k+\rho-3}.\]

This gives the desired conclusion about the structure of $A^{red}_k$.

We now examine the support of the summand $\mathcal{T}_{-2V_k}$ of $A^+_k$.

\begin{claim}\label{claim:2Vk}
    The lowest grading in $\mathcal{T}_{-2V_k}$ is of the form $-2V_k = \min(0,k+\rho-\eta\pm 1)$, where $\eta = 1 + (k+\rho \mod 2)$.
\end{claim}
\begin{proof}[Proof of Claim \ref{claim:2Vk}]
    If $k+\rho-3 \geq 0$, then $-2V_k=0$ by the proof of Claim \ref{claim:Aredmaxsupport}. Since $k+\rho-\eta\pm 1 \geq k+\rho-3$, we have $\min(0,k+\rho-\eta\pm 1) = 0=-2V_k$ as desired.
    
    If $k+\rho-3 < 0$, then $A^+_k$ is supported in degrees at least $k+\rho-3$ by the proof of Claim \ref{claim:Aredmaxsupport}. Therefore, $-2V_k \geq k+\rho-3$.

    In grading $k+\rho-1$, the exact triangle (\ref{triangle:injection}) and Lemma \ref{lemma:basicFU} give an injection
    \begin{equation}\label{diag:injection}
    \begin{tikzcd}
    \mathcal{T}|_{k+\rho-1}  \arrow[r, hook] & \mathcal{T}_{-2V_k}|_{k+\rho-1} \subset A^+_k|_{k+\rho-1}
    \end{tikzcd}.
    \end{equation}

    It follows that $\mathcal{T}_{-2V_k}$ is non zero in each even grading greater than or equal to $k+\rho-1$. Therefore, $k+\rho-3 \leq -2V_k \leq k+\rho-\eta+1$ and we obtain $-2V_k = k+\rho-\eta\pm 1$. Further, if $k+\rho-3 \geq -2$, then $k+\rho-\eta+ 1 > 0 \geq -2V_k$ and thus $-2V_k =k+\rho-\eta-1=\min(0,k+\rho-\eta- 1)$. If $k+\rho-3 < -2$, then both $k+\rho-\eta+ 1=\min(0,k+\rho-\eta+ 1)$ and $k+\rho-\eta- 1=\min(0,k+\rho-\eta- 1)$ are possible values for $-2V_k$.
    \end{proof}

This concludes the proof of the lemma.
\end{proof} 


\subsection{Property \spl{} for thickness-two knots}\label{subsec:th2spliff} 

Lemma \ref{lemma:structth2} says that $A^{red}_k$ is of the form
\[\mathcal{T}_{2k}|_{\leq k+\rho-2 + \epsilon} \oplus \left(\frac{\F_{k+\rho-1}[U]}{U^2}\right)^r \oplus\ \F^a_{k+\rho-1} \oplus \F_{k+\rho-2}^b \oplus \F^c_{k+\rho-3},\]
where $\epsilon = \pm 1$. We now examine each possibility for $\epsilon, a,b,c$ and verify whether $A^+_k$ has property \spl{}. Note that we can ignore $r$ since $1 \in \F[U]/U^2$ is not killed by the $U$-action. If both $a$ and $c$ are non-zero, then $A^+_k$ does not have property \spl{}. We may thus assume that at least one of $a$ or $c$ is zero.

First, suppose $k+\rho$ is odd. Elements of odd degree may only appear in $\F_{k+\rho-2}^b$, so we are interested only in the values of $\epsilon, a$ and $c$.
\begin{itemize}
    \item If $a=c=0$, then $A^+_k$ has property \spl{}.
    \item If $a\neq 0$ and $c=0$, then $A^+_k$ does not have property \spl{} if and only if $\mathcal{T}_{2k}|_{\leq k+\rho-2+ \epsilon}$ is generated by a unique element of degree $k+\rho-3$. This happens if and only if $k = \rho-3$ and $\epsilon = -1$.
    \item If $a=0$ and $c\neq0$, then $A^+_k$ does not have property \spl{} if and only if $\mathcal{T}_{2k}|_{\leq k+\rho-2+ \epsilon}$ is generated by a unique element of degree $k+\rho-1$. This happens if and only if $k = \rho-1$ and $\epsilon = +1$.
\end{itemize}

Suppose now that $k+\rho$ is even. Elements of odd degree may only appear in one of $\F_{k+\rho-1}^a$ or $\F_{k+\rho-3}^c$, so we are interested only in the values of $\epsilon$ and $b$.
\begin{itemize}
    \item If $b=0$, then $A^+_k$ has property \spl{}.
    \item If $b \neq 0$ and $\epsilon = +1$, then $A^+_k$ has property \spl{} because the even elements of $A^{red}_k$ that are not in $\F_{k+\rho-2}^b$ must appear in $\mathcal{T}_{2k}|_{\leq k+\rho-2}$, which has non-zero $U$-action or is supported in degree $k+\rho-2$ if it is non-trivial.
    \item If $b \neq 0$ and $\epsilon = -1$, then $A^+_k$ fails to have property \spl{} if and only if $\mathcal{T}_{2k}|_{\leq k+\rho-4}$ is supported only in degree $k+\rho-4$. This happens if and only if $k = \rho-4$.
\end{itemize}

This is summarized in Table \ref{table:spliFf-th2}.

\begin{table}[ht!]
\begin{tabular}{|l||l|l|l|l|l|}
\hline
\multirow{2}{*}{$(a,c)$} & \multicolumn{2}{c}{$k+\rho$ odd} & \multicolumn{3}{|c|}{$k+\rho$ even} \\\cline{2-6}
                     & $\epsilon=+1$   &$\epsilon=-1$ & $b=0$              & $\epsilon=+1, b\neq 0$& $\epsilon=-1,b\neq0$ \\\hline
$(1,0)$ & yes                   & yes iff $k\neq\rho-3$                   &yes              & yes                   & yes iff $k\neq\rho-4$ \\
$(0,1)$ & yes iff $k\neq\rho-1$                 & yes                    & yes              & yes                   & yes iff $k\neq\rho-4$ \\
$(0,0)$ & yes                   & yes                   &yes              & yes                   & yes iff $k\neq\rho-4$ \\
$(1,1)$ & no                    & no                    &no               & no                    & no    \\\hline
\end{tabular}
\caption{Structure of $A^+_k$ and satisfaction of property \spl{} for knots of thickness two}
\label{table:spliFf-th2}
\end{table}

Using the maps (\ref{diag:surjection}) and (\ref{diag:injection}), we can guarantee that $A^+_k$ has property \spl{} given certain conditions on $\HFK(K)$.

\begin{lemma}\label{lemma:rho-3rho-4}
Let $K$ be a knot of thickness two. Let $\rho$ be an integer such that for all $s$, the group $\HFK_d(K,s)$ is non-zero only for gradings $d \in \{s+\rho, s+\rho-1, s+\rho-2\}$. Suppose $k \geq 0$.
\begin{enumerate}[(i)]
    \item If $k+\rho$ is odd, then $A^+_k$ has property \spl{} if  $\HFK_{k+\rho}(K,k)=0$.
    \item If $k+\rho$ is odd and $k \neq \rho-3$, then $A^+_k$ has property \spl{} if  $\HFK_{k+\rho-2}(K,k)=0$.
    \item If $k+\rho$ is even and $k \neq \rho-4$, then $A^+_k$ has property \spl{} if at least one of the groups $\HFK_{k+\rho}(K,k)$ or $\HFK_{k+\rho-2}(K,k)$ is trivial.
\end{enumerate}
\end{lemma}
\begin{proof}
     If $\HFK_{k+\rho}(K,k)=0$, then $C_{\{i \leq -1 \wedge j \leq k - 1\}}$ does not contain elements of degree $k+\rho-2$, so neither does its homology group. Hence, the injection (\ref{diag:injection}) is an isomorphism and $A^{red}_k$ is of the form
     \[\mathcal{T}_{2k}|_{\leq k+\rho-3} \oplus \F_{k+\rho-2}^b  \oplus \F_{k+\rho-3}^c.\]
     This corresponds to columns 3 to 6 of the $(0,1)$ row of Table \ref{table:spliFf-th2}. Hence, $A^+_k$ fails to have property \spl{} only if $k+\rho$ is even and $k=\rho-4$.

    If $\HFK_{k+\rho-2}(K,k)=0$, then $C_{\{i \geq 0 \wedge j \geq k\}}$ does not contain elements of degree $k+\rho-2$, so neither does its homology group. Hence, the surjection (\ref{diag:surjection}) is an isomorphim and $A^{red}_k$ is of the form
     \[\mathcal{T}_{2k}|_{\leq k+\rho-2 \pm 1} \oplus \F_{k+\rho-1}^a \oplus \F_{k+\rho-2}^b.\]
     This corresponds to the $(1,0)$ row of Table \ref{table:spliFf-th2}. Hence, $A^+_k$ fails to have property \spl{} only if $k=\rho-3$ when $k+\rho$ is odd, and only if $k=\rho-4$ when $k+\rho$ is even.
\end{proof}

We now turn to the proof of Proposition \ref{prop:th2hfkspliff}.

\begin{propth2hfkspliff}
Let $K$ be a knot of thickness two. Let $\rho$ be an integer such that for all $s$, the knot Floer homology group $\HFK_d(K,s)$ is non-zero only for gradings $d \in \{s+\rho, s+\rho-1, s+\rho-2\}$.

    Suppose $\rho \in \{0,1,2\}$. If for each $k \geq 0$, at least one of the groups $\HFK_{k+\rho}(K,k)$ or $\HFK_{k+\rho-2}(K,k)$ is trivial, then $K$ and its mirror both satisfy property \spl{}. Therefore, $K$ admits only finitely many non-integral non-characterizing Dehn surgeries.
 \end{propth2hfkspliff}   

We first show a slightly more general statement for $K$.

\begin{lemma}\label{lemma:rho012-K}
Let $K$ be a knot of thickness two. Let $\rho$ be an integer such that for all $s$, the knot Floer homology group $\HFK_d(K,s)$ is non-zero only for gradings $d \in \{s+\rho, s+\rho-1, s+\rho-2\}$.

    Suppose $\rho \leq 2$. If for each $k \geq 0$, at least one of the groups $\HFK_{k+\rho}(K,k)$ or $\HFK_{k+\rho-2}(K,k)$ is trivial, then $K$ has property \spl{}.
\end{lemma}
\begin{proof}
We have $\rho-3, \rho-4 < 0$, so $k \neq \rho-3, \rho-4$ for all $k \geq 0$. By Lemma \ref{lemma:rho-3rho-4}, $A^+_k$ has property \spl{} for all $k \geq 0$. Hence, $K$ has property \spl{}.
\end{proof}

\begin{proof}[Proof of Proposition \ref{prop:th2hfkspliff}]
The statement for $K$ follows from Lemma \ref{lemma:rho012-K}, so we need to show that $mK$ also has property \spl{}. Recall the symmetry properties of knot Floer homology \cite{osHF}
\begin{itemize}
    \item[(S1)] $\HFK_d(K,s)\cong \HFK_{-d}(mK,-s)$ and
    \item[(S2)] $\HFK_d(K,s)\cong \HFK_{d-2s}(K,-s)$.
\end{itemize}
Let $\rho_{m}$ denote the integer such that for all $s$, the group $\HFK(mK,s)$ is non-zero only for gradings $d \in \{s+\rho_m, s+\rho_m-1, s+\rho_m-2\}$. By (S1), we have $\rho_m = 2-\rho \in \{0,1,2\}$. Further, by (S1) and (S2), we have isomorphisms $\HFK_{k+\rho_m-2}(mK,k)\cong \HFK_{k+\rho}(K,k)$ and $\HFK_{k+\rho_m}(mK,k)\cong \HFK_{k+\rho-2}(K,k)$. 

Therefore, $mK$ satisfies the hypotheses of Lemma \ref{lemma:rho012-K} and thus has property \spl{}. It follows from Theorem \ref{theor:mccspliff} that any knot of thickness two that satisfies the assumptions of Proposition \ref{prop:th2hfkspliff} verifies Conjecture \ref{conj:mcc}. 
\end{proof}

\subsection{Computations for thickness-two knots}\label{subsec:th2compute}
To verify if Conjecture \ref{conj:mcc} holds for a knot of thickness two, we need to compute the structure of its modules $A^+_k$ for all $k \geq 0$ that do not satisfy the conditions of Lemma \ref{lemma:rho-3rho-4}. To achieve this, one may first compute $CFK^\infty(K)$, following the approach used for thickness-one knots. However, two main issues arise when dealing with knots of thickness greater than one. 

\subsubsection{Computing lifts}
First, it is difficult in general to find a lift of $CFK_{\F[u,v]}(K)/(uv)$, in the sense of Section \ref{section:overview}. We can still exploit the computational effectiveness of solving linear systems, as was done in the case of thickness-one knots, to reduce the number of possibilities for diagonal arrows. Recall that we encode the unknown differential map acting on the underlying $\FU$-module $C$ of $CFK_{\F[u,v]}(K)/(uv)$ as a matrix $d_{var}=H+V+D_{var}$. As in Algorithm \ref{algo:cfk1}, we obtain a system of equations $E$ by setting $d_{var}^2=0$, but it may contain non-linear equations if Proposition \ref{prop:linear} is not satisfied. By considering the maximal subsystem of linear equations $E'$ of $E$, we obtain a matrix equation $Aa=b$ with an initial solution $a=a_0$ and whose set of solutions is $a_0 + \ker A$. If $E' \neq E$, we need to determine which elements of $a_0 + \ker A$ are solutions of the full system $E$. Indexing the elements of $a_0 + \ker A$ by $a_l, l=0, \ldots, 2^{\dim \ker A}-1$, we obtain maps $d_l = H+V+D_{var}(a_l) = d_{var}(a_l)$. If $d_l^2=0$, then $a_l$ is a solution to $E$ and the differential complex $\mathcal{C}_l = (C, d_l)$ is a lift of $CFK_{\F[u,v]}(K)/(uv)$. Note that this approach is computationally manageable only when the dimension of $\ker A$ is relatively small, or when Proposition \ref{prop:linear} is satisfied, in which case the set of lifts is $\{\mathcal{C}_l \,|\, a_l \in a_0 + \ker A$\}.

\subsubsection{Equivalence of lifts}
Second, the computed lifts may not be filtered chain homotopy equivalent to one another. For knots with up to 16 crossings, Hanselman showed that $CFK_{\F[u,v]}(K)$ splits as in Theorem \ref{theor:popovic} \cites[Corollary 12.6]{hanselman}{hanselman_website}; therefore, any lift of $CFK_{\F[u,v]}(K)/(uv)$ is a genuine representative of the full knot Floer complex. It then suffices to verify property \spl{} for the modules $A^+_k$ of any lift obtained by the method described previously. This is done by using Lemma \ref{lemma:rho-3rho-4} and adapting the method described in Section \ref{section:Akimpl} to thickness-two knots, according to Lemma \ref{lemma:structth2}.

For knots with at least 17 crossings, we may not have such an equivalence between lifts. However, for our application at hand, we are interested only in the modules $A^+_k$, which are the homology groups of quotients of $CFK^\infty(K)$. Our strategy thus consists in computing all possible lifts of $CFK_{\F[u,v]}(K)/(uv)$ by considering each element in the set $a_0 + \ker A$. We then check that the modules $A^+_k$ of each of these lifts -- which may belong to different homotopy equivalence classes -- verify property \spl{}. Table \ref{table:17th2spliff} summarizes the results of our computation for knots with 17 crossings, carried out for knots and their mirrors whose maximal linear subsystems have kernel of dimension at most 12. For each dimension, the table indicates the number of knots with both the knot and its mirror satisfying property \spl{}, and the number with either the knot or its mirror not satisfying property \spl{}.

In all computed cases, every lift for a given knot gave the same outcome for the presence or absence of property \spl{}. We did not verify whether the different lifts were chain homotopy equivalent, but inspection of the relevant groups $A^+_k$ for some of these knots revealed that they have the same graded $\F[U]$-module structure across all lifts. Therefore, for these cases, the quotient complex $CFK_{\F[u,v]}(K)/(uv)$ appears to fully determine the homology groups $A^+_k$. 
This observation suggests that this, or the stronger statement of Theorem \ref{theor:introcfk}, could hold not only for thickness-one knots, but also for thicker knots whose quotient complex $CFK_{\F[u,v]}(K)/(uv)$ is sufficiently simple. It remains unclear what precise algebraic conditions would convey this simplicity, or whether such conditions exist at all beyond thickness one.

\begin{table}[ht!]
	\centering
	\begin{tabular}{|c|cc|}
		\hline
		$\dim \ker A$  & $K$ and $mK$ \spl{} & $K$ or $mK$ non-\spl{} \\ \hline
		0 & 498 & 0  \\
		2 & 174 & 6  \\
		4 & 155 & 20 \\
		6 & 153 & 21 \\ 
		8 & 117 & 39 \\ 
		10 & 95 & 26 \\ 
		12 & 135 & 31 \\ \hline
	\end{tabular}
	\caption{Thickness-two knots with 17 crossings and property \spl{}} \label{table:17th2spliff}
\end{table}


\bibliographystyle{amsalpha} 
\bibliography{bib}

\end{document}